\documentclass[11pt,leqno]{amsart}
\usepackage{amsmath,amssymb,array,latexsym}

\newtheorem{thm}{Theorem}[section]
\newtheorem{lem}[thm]{Lemma}
\newtheorem{prop}[thm]{Proposition}
\newtheorem{cor}[thm]{Corollary}

\theoremstyle{definition}
\newtheorem*{defin}{Definition}
\newtheorem*{remark}{Remark}

\newcommand{\N}{\mathbb{N}}
\newcommand{\R}{\mathbb{R}}

\newcommand{\cA}{\mathcal{A}}

\def\trivert#1{|\!|\!| #1|\!|\!|}
\def\Btrivert#1{\Big|\!\Big|\!\Big| #1\Big|\!\Big|\!\Big|}

\newcommand{\vp}{\varepsilon}
\newcommand{\supp}{{\rm supp}}

\newcommand{\coo}{\mathrm{c}_{00}}
\newcommand{\co}{\mathrm{c}_0}

\newcommand{\cof}{{\rm cof}}

\newcommand{\bb}{\overline{b}}

\newcommand{\xb}{\overline{x}}

\def\hangbox to #1 #2{\vskip1pt\hangindent #1\noindent \hbox to #1{#2}$\!\!$}

\newcommand{\norm}[1]{\lVert#1\rVert}
\newcommand{\abs}[1]{\lvert#1\rvert}

\newcommand{\bignorm}[1]{\bigl\lVert#1\bigr\rVert}

\newcommand{\Bignorm}[1]{\Bigl\lVert#1\Bigr\rVert}
\newcommand{\biggnorm}[1]{\biggl\lVert#1\biggr\rVert}

\newcommand{\vpb}{\overline\varepsilon}
\newcommand{\deltab}{\overline\delta}
\newcommand{\etab}{\overline\eta}

\newcommand{\Et}{\tilde E}
\newcommand{\xt}{\tilde x}
\newcommand{\yt}{\tilde y}

\newcommand{\wt}{\tilde w}

\newcommand{\Qt}{\tilde Q}

\newcommand{\Gt}{\tilde G}

\newcommand{\Yt}{\tilde Y}

\newcommand{\Zs}{Z^{(*)}}
\newcommand{\Us}{U^{(*)}}
\newcommand{\Vs}{V^{(*)}}

\allowdisplaybreaks[4]

\author{E. Odell, Th. Schlumprecht and A.~Zs\'ak}
\thanks{Research of the first two authors was supported by the
National Science Foundation}

\title[Asymptotic $\ell_p$ spaces]{On the structure of asymptotic
  $\ell_p$ spaces}

\subjclass[2000]{Primary 46B20; Secondary 46B03, 46B10}
\begin{document}

\date{20/2/2006}

\begin{abstract}
  We prove that if $X$ is a separable, reflexive space which is
  asymptotic $\ell_p$ for some $1\le p\le \infty$, then $X$ embeds into
  a reflexive space $Z$ having an asymptotic $\ell_p$
  finite-dimensional decomposition. This result leads to an intrinsic
  characterization of subspaces of  spaces with an  asymptotic
  $\ell_p$ FDD. More general results of this type are also obtained.
\end{abstract}
\maketitle
\section{Introduction}\label{S:1}

Let $X$ be a separable Banach space with a finite-dimensional decomposition
(FDD), $(E_n)$. Let $1\le p\le \infty$. $X$ is {\em asymptotic
  $\ell_p$  with respect to $(E_n)$} \cite{MT} if there exists
$C<\infty$ so that for all $n$ and all block sequences $(x_i)_{i=1}^n$
of $(E_i)_{i=n}^\infty$
$$ \frac1C\Big(\sum_{i=1}^n \|x_i\|^p\Big)^{1/p}\le
\Big\|\sum_{i=1}^n x_i\Big\|\le C\Big(\sum_{i=1}^n \|x_i\|^p\Big)^{1/p}$$
(if $p=\infty$ we use the $c_0$-norm $\max \|x_i\|$).

A {\em coordinate-free version} of this notion is as follows
\cite{MMT}. Let $X$ be an arbitrary Banach space, and let $\cof(X)$
denote the set of all closed  subspaces of $X$ having finite
codimension.
We say $X$ is {\em asymptotic $\ell_p$ }if there exists $C<\infty$ so that
\begin{align}\label{E:1.1}
  \forall n\in\N\quad&\exists Y_1\in \cof(X)\ \forall y_1\in S_{Y_1}
  \text{ (unit sphere of $Y_1$)}\\
  &\exists Y_2\in \cof(X)\ \forall y_2\in S_{Y_2} \notag \\
  &\vdots \notag\\
  &\exists Y_n\in \cof(X)\ \forall y_n\in S_{Y_n} \notag\\
  &\qquad (y_i)_{i=1}^n \text{ is $C$-equivalent to the unit vector
    basis of $\ell_p^n$.} \notag
\end{align}

We prove (Corollary~\ref{C:4.6}) that if $X$ is separable, reflexive
and asymptotic $\ell_p$, then there exists a reflexive space  $Z$
with an FDD $(E_n)$ so that $Z$ is asymptotic $\ell_p$ with respect
to $(E_n)$ and $X$ embeds isomorphically into $Z$.  This is deduced
from a more general result (Corollary~\ref{C:4.7}) which considers
separable, reflexive spaces $X$ that satisfy the following for some
$1\le q\le p \le\infty$ and $C<\infty$.
\begin{align}\label{E:1.2} 
  \forall n\in\N\quad&\exists Y_1\in \cof(X)\ \forall y_1\in S_{Y_1} \\
  &\exists Y_2\in \cof(X)\ \forall y_2\in S_{Y_2} \notag \\
  &\vdots \notag\\
  &\exists Y_n\in \cof(X)\ \forall y_n\in S_{Y_n} \notag\\
  &\quad
  \frac1C \Big(\sum_{i=1}^n |a_i|^p\Big)^{1/p}\le
  \Big\|\sum_{i=1}^n a_i y_i\Big\|
  \le C \Big(\sum_{i=1}^n |a_i|^q\Big)^{1/q}
  \text{ whenever }(a_i)_{i=1}^n\subset\R\ .\notag
\end{align}
We characterize such spaces as those that embed into reflexive spaces
$Z$ with an FDD $(E_i)$ {\em satisfying asymptotic
  $(\ell_p,\ell_q)$-estimates}. This means that for some $C<\infty$,
for all $n\in\N$ and all block sequences $(x_i)_{i=1}^n$ of
$(E_i)_{i=n}^\infty$ $$ \frac1C\Big(\sum_{i=1}^n
\|x_i\|^p\Big)^{1/p}\le\Big\|\sum_{i=1}^n x_i\Big\|\le
C\Big(\sum_{i=1}^n \|x_i\|^q\Big)^{1/q}.$$ We also show that this is
equivalent to $X$ being a quotient of such a space $Z$.

To accomplish that we develop a more general machinery concerning
asymptotic $U$-upper and $V$-lower estimates, where $U$ and $V$ are
certain spaces with subsymmetric bases. Theorem~\ref{T:4.6} is
concerned with this general setting, and Corollary~\ref{C:4.7} will
then be obtained as a special case.

In an earlier paper \cite{OS1} analogous  results were obtained
characterizing when a reflexive space embeds into
$\big(\oplus_{i=1}^\infty E_i\big)_{\ell_p}$, the $\ell_p$-sum of some
sequence $E=(E_i)$ of finite-dimensional spaces. In this paper the
role of $\big(\oplus_{i=1}^\infty E_i\big)_{\ell_p}$ is played by a
space $Z_V(E)$, where $V$ is a space with an unconditional basis,
especially a convexified Tsirelson space. Section \ref{S:2} contains
the precise definition of $Z_V(E)$ and some structural results about
it. In  section \ref{S:3} we present several   embedding theorems
(Theorem \ref{T:3.1}, Corollary \ref{C:3.2} and Theorem \ref{T:3.1b})
which characterize when a reflexive space embeds into some $Z_V(E)$
and is a quotient of some $Y_V(F)$. Section \ref{S:4} contains our
main results.

Before proceeding we  first introduce some
definitions and  notation.

Let $Z$ be a Banach space with an FDD $E=(E_n)$.  For $n\in\N$ we denote
the $n$-th {\em coordinate projection } by $P^E_n$,
i.e. $P_n^E:Z\to E_n,\quad \sum z_i\mapsto z_n.$ For  finite
$A\subset \N$ we put $P^E_A=\sum_{n\in A} P_n^E$. The  {\em projection
  constant  } of $(E_n)$ (in $Z$) is defined by
$$K=K(E,Z)=\sup_{m\le n}\|P_{[m,n]}^E\|.$$
Recall that $K$ is always finite and, as in the case of bases, we call
{\em  $(E_n)$ bimonotone (in $Z$)} if $K=1$.  By passing to the
equivalent norm
$$\trivert{\cdot}:Z\to\R,\quad z\mapsto\sup_{m\le n} \|P^E_{[m,n]}(z)\|,$$
we can always renorm $Z$ so that $K=1$.

For a sequence $(E_i)$ of finite-dimensional spaces
we define the vector space
$$\coo(\oplus_{i=1}^\infty E_i)=\big\{(z_i):z_i\in E_i
\text{ for $i\in\N$, and $\{i\in\N: z_i\not=0\}$ is finite}\big\},$$
which is dense in each Banach space for which $(E_i)$ is an FDD. For
$A\subset \N$ we denote  by $\oplus_{i\in A} E_i$ the linear subspace
of $\coo(\oplus E_i)$ generated by the elements
of $\bigcup_{i\in A}E_i$. As usual we denote the vector space of sequences
in $\R$ which are eventually zero by $c_{00}$ and its  unit vector
basis by $(e_i)$.
We  sometimes will consider for the same sequence $(E_i)$ of
finite-dimensional spaces different norms on $\coo(\oplus
E_i)$. In order to avoid confusion we will therefore often index the
norm by the Banach space whose norm we are using, i.e. $\|\cdot\|_Z$
denotes the norm  of the Banach space $Z$.

If $Z$ has an FDD $(E_i)$, the vector space $\coo(\oplus_{i=1}^\infty
E^*_i)$, where $E^*_i$ is the dual space of $E_i$ for $i\in\N$, is a
$w^*$-dense subspace of $Z^*$. We denote the norm closure of
$\coo(\oplus_{i=1}^\infty E^*_i)$ in $Z^*$ by $\Zs$. $\Zs$ is
$w^*$-dense in $Z^*$, the unit ball $B_{\Zs}$ norms $Z$, and $(E_i^*)$
is an FDD of $\Zs$ having a projection constant not exceeding
$K(E,Z)$. If $K(E,Z)=1$ then $B_{Z^{(*)}}$ is 1-norming and
$Z^{(*)(*)}=Z$.

For $z\in \coo(\oplus E_i)$ we define the $E$-{\em support of $z$} by
\[
\supp_E(z)=\{i\in\N:  P^E_i(z)\not=0\}.
\]
A sequence $(z_j)$ (finite or infinite) of non-zero vectors in
$\coo(\oplus E_i)$ is called  a {\em block sequence of $(E_i)$} if
$$\max\supp_E(z_n)<  \min\supp_E(z_{n+1})\qquad\text{whenever $n\in\N$ (or
  $n<$length$(z_j)$)},$$ and it is called a  {\em skipped block
  sequence of $(E_i)$} if
$$\max\supp_E(z_n)<  \min\supp_E(z_{n+1})-1\qquad\text{whenever $n\in\N$ (or
  $n<$length$(z_j)$)}.$$
Let $\deltab=(\delta_n)\subset (0,1)$ with $\delta_n\downarrow 0$.
A (finite or infinite) sequence $(z_j)\subset S_Z=\{z\in Z:
\|z\|=1\}$ is called a $\deltab$-{\em block sequence of $(E_n)$}
or a $\deltab$-{\em skipped block sequence of $(E_n)$}
if
there are $k_1<k_2<\ldots$ in $\N$ so that ($k_0=1$)
\[
\| z_n - P^E_{[k_{n-1},k_n)}(z_n)\|<\delta_n,\text{ or }\| z_n -
P^E_{(k_{n-1},k_n)}(z_n)\|<\delta_n, \text{ respectively,}
\]
for all $n\in\N$  (or $n\le$length$(z_j)$).
Of course one could generalize the notion of $\deltab$- block and
$\deltab$-skipped block sequences to more general sequences, but
we prefer to introduce this notion only for normalized sequences.

\begin{remark}
  If $(F_i)$ is a blocking of $(E_i)$ and if $(x_i)$ is a
  $\deltab$-skipped block sequence of $(F_i)$, then $(x_i)$ is not
  necessarily a $\deltab$-skipped block sequence of $(E_i)$ (since  in
  the
  definition of skipped block sequence we skip exactly one
  coordinate). Nevertheless it is clear that $(x_i)$ is a
  $2K\deltab$-skipped block sequence of $(E_i)$, where $K$ is the
  projection constant of $(E_i)$ in $Z$.
\end{remark}

A sequence of finite-dimensional spaces $(G_n)$ is  called
  a {\em blocking of $(E_n)$} if
 there are
$0=k_0<k_1<k_2<\ldots$ in $\N$ so that
$G_n=\oplus_{i=k_{n-1}+1}^{k_n} E_i$ for $n=1,2,\ldots$.

\begin{defin}
  \label{D:1.1}
  For two normalized basic sequences $(e_i)$ and $(f_i)$ we say that
  {\em $(f_i)$ $C$-dominates  $(e_i)$} or that \emph{$(e_i)$   is
  $C$-dominated by $(f_i)$}, where   $C \ge 1$, if for all $(a_i)\in
  \coo$
  $$ \Bigl\|\sum a_ie_i\Bigr\|\le C  \Bigl\|\sum a_if_i\Bigr\|.$$
  
  We say that {\em $(f_i)$ dominates  $(e_i)$} or that \emph{$(e_i)$
  is dominated by $(f_i)$} if, for some $C\ge 1$,  $(f_i)$
  $C$-dominates  $(e_i)$.

  Let $V$ be a Banach space with a
  1-unconditional and normalized basis $(v_i)$ and let $1\le C<\infty$.
  We say that
  an FDD $(E_n)$ of a Banach space $Z$ {\em satisfies  $C$-$V$-lower
    estimates (in $Z$)}   if  all normalized
  block sequences of $(E_n)$ in $Z$   $C$-dominate   $(v_i)$,
  and   $(E_n)$  {\em satisfies  $C$-$V$-upper estimates (in
    $Z$)} if $(v_i)$ $C$-dominates all normalized block sequences of
  $(E_i)$ in $Z$.
  If  $U$ is another space with a normalized and 1-unconditional basis
  $(u_i)$, we say that $(E_n)$ {\em  satisfies  $C$-$(V,U)$-estimates
    (in $Z$)} if it   satisfies $C$-$V$-lower estimates and
  $C$-$U$-upper estimates.

  We say that $(E_n)$ {\em  satisfies $V$-lower estimates (in $Z$)},
  {\em  $U$-upper estimates (in $Z$)} or
  {\em $(V,U)$-estimates  (in $Z$)} if
  there is a constant $C$ so that $(E_n)$   satisfies $C$-$V$-lower estimates (in $Z$),
  $C$-$U$-upper estimates (in $Z$) or
  $C$-$(V,U)$-estimates  (in $Z$), respectively.
\end{defin}
\begin{remark}
  It is easy to show that if every normalized block sequence of
  $(E_i)$ in $Z$ dominates $(v_i)$, then $(E_i)$ satisfies
  $V$-lower estimates in $Z$. A similar remark holds for
  $U$-upper estimates.
\end{remark}

 We define  for $\ell\in\N$
\begin{align*}
T_\ell&=\big\{(n_1,n_2,\ldots,n_\ell): n_1<n_2<\ldots<n_\ell
 \text{ are in }\N\big\}
\intertext{and}
T_\infty&=
\bigcup_{\ell=1}^\infty T_\ell .
\end{align*}
If $\alpha=(m_1,m_2,\ldots,m_\ell)\in T_\ell$, we call $\ell$ the {\em length of $\alpha$}
 and denote it by $|\alpha|$, and $\beta=(n_1,n_2,\ldots,n_k)\in T_\infty$ is called
  an {\em extension of $\alpha$}, or $\alpha$ is called {\em a restriction of }
$\beta$, if $k \ge \ell$ and $n_i=m_i$ for $i=1,2,\ldots,\ell$.
 We then write $\alpha\le \beta$ and with this order $(T_\infty,\le)$ is a tree.

A set $S\subset T_\infty$ is called {\em well-founded} if
 it is closed under taking restrictions  and if it does not contain
 any infinite chain with respect to $\le$. Note that this means
that the set $\max(S)$ of maximal elements of $S$ is not empty
(provided $S\not=\emptyset$) and that
 $$S=\big\{ (n_1,n_2,\ldots,n_k):\exists \ell\ge k\,\exists n_k<n_{k+1}<\ldots <n_\ell \quad
 (n_1,\ldots,n_\ell)\in\max(S)\big\}.$$

In this work  {\em trees } in a Banach space $X$ are families  in $X$
indexed by
$T_\infty$, and thus they are countable infinitely branching trees
 of countably infinite length.
 In Section \ref{S:4} we will consider families in $X$ indexed by $T_\ell$ for some
 $\ell\in\N$, and we refer to them as {\em trees of length $\ell$}.

 We introduce the following notation only for trees of infinite length, but note
 that they can be similarly defined for trees of finite length.

 For a  tree
 $(x_\alpha)_{\alpha\in T_\infty}$ in a  Banach space $X$,
 and $\alpha=(n_1,n_2,\ldots, n_\ell)\in T_\infty\cup\{\emptyset\}$
  we call the sequences of the form $(x_{(\alpha,n)})_{n>n_\ell}$
 {\em nodes of $(x_\alpha)_{\alpha\in  T_\infty}$}.
 The sequences $(y_n)$ with $y_i=x_{(n_1,n_2,\ldots,n_i)}$ for
 $i\in\N$ and for some strictly increasing sequence $(n_i)\subset\N$
 are called {\em branches of} $(x_\alpha)_{\alpha\in T_\infty}$. Thus
 branches of a tree $(x_\alpha)_{\alpha\in T_\infty}$ are sequences of
 the form $(x_{\alpha_n})$, where $(\alpha_n)$ is an increasing (with
 respect to extension) sequence in $T_\infty$ with $|\alpha_n|=n$
 for all $n\in\N$.

 If $(x_\alpha)_{\alpha\in T_\infty}$  is a tree in $X$ and if $T'\subset T_\infty$  is closed under
 taking restrictions so that for each $\alpha\in T'\cup\{\emptyset\}$
 infinitely many direct successors of $\alpha$ are also in $T'$, then we
call $(x_\alpha)_{\alpha\in T'}$ a {\em full subtree of }
$(x_\alpha)_{\alpha\in T_\infty}$. Note that
  $(x_\alpha)_{\alpha\in T'}$ could then be relabeled to a family indexed by
$T_\infty$ and note that the branches of $(x_\alpha)_{\alpha\in T'}$
are branches of
 $(x_\alpha)_{\alpha\in T_\infty}$
 and that the nodes of $(x_\alpha)_{\alpha\in T'}$ are subsequences of  certain nodes of
 $(x_\alpha)_{\alpha\in T_\infty}$.

 We call a tree $(x_\alpha)_{\alpha\in T_\infty}$  in a  Banach space $X$
 {\em normalized} if $\|x_\alpha\|=1$ for all  $\alpha\in T_\infty$ and
{\em weakly null} if every node
 is weakly null.
    If  $(x_\alpha)_{\alpha\in T_\infty}$ is a tree in a Banach space
 $Z$
which
has an FDD $(E_n)$, then  we call it a
 {\em  block tree  of } $(E_n)$ if every node is a
 block sequence of $(E_n)$.

We shall need a coordinate-free version of lower and upper estimates.

\begin{defin}
  \label{D:1.2}
  Let $V$ be a Banach space with a 1-unconditional and normalized
 basis $(v_i)$ and let $1\le C<\infty$. We say that a Banach space $X$
 {\em satisfies $C$-$V$-lower tree estimates }if every normalized
 weakly null tree $(x_\alpha)_{\alpha\in T_\infty}$ in $X$ has a
 branch $(y_i)$ which $C$-dominates  the basis $(v_i)$. Of course we
 defined domination for basic sequences only but since every
 normalized, weakly null tree  in $X$ admits a full subtree with all
 branches $2$-basic, say, this does not constitute a problem.

 We say that $X$  {\em satisfies   $C$-$V$-upper tree estimates }if
 every normalized weakly null tree $(x_\alpha)_{\alpha\in T_\infty}$
 in $X$ has a branch $(y_i)$  which is  $C$-dominated  by $(v_i)$.
 
 If $U$ is a second space with a 1-unconditional and normalized basis
 $(u_i)$,
 we say that  $X$
 {\em satisfies $C$-$(V,U)$-tree estimates} if
 it satisfies $C$-$V$-lower and $C$-$U$-upper tree estimates.
 
 We say that $X$ {\em satisfies $V$-lower tree, $U$-upper tree} or
 {\em $(V,U)$-tree estimates} if for some $1\le C<\infty$ $X$ satisfies
 $C$-$V$-lower tree, $C$-$U$-upper tree or $C$-$(V,U)$-tree estimates,
 respectively.
\end{defin}

\begin{prop}
  \label{tree-est}
  Let $U$ and $V$ be Banach spaces with normalized, 1-unconditional
  bases $(u_i)$ and $(v_i)$, respectively.  For an infinite subset
  $N\subset \N$ we let $U^{(N)}$ and $V^{(N)}$  be the closed
  subspaces spanned by $(u_{i})_{i\in N}$ and $(v_{i})_{i\in N}$,
  respectively.

  If $C\ge 1$ and a Banach space $X$ satisfies  $C$-$(V,U)$-tree
  estimates, then it also satisfies $C$-$(V^{(N)},U^{(N)})$-tree
  estimates (with respect to the 1-unconditional bases $(v_{i})_{i\in
    N}$ and $(u_{i})_{i\in N}$) for any infinite $N$.
\end{prop}
\begin{proof}
  Let $n_1<n_2<n_3<\ldots $ be such that $N=\{n_1,n_2,\ldots\}$, and
  let $(x_\alpha)_{\alpha\in T_\infty}$ be a normalized, weakly null
  tree in $X$. Let $(z_n)$ be any  weakly null sequence in $S_X$
  (e.g.~the top node of $(x_\alpha)_{\alpha\in T_\infty}$).
  
  We now consider
  the following tree  $(\xt_\alpha)_{\alpha\in T_\infty}$
  which, up to finitely many elements of each node, is an expansion
  of  $(x_\alpha)_{\alpha\in T_\infty}$: for
  $\alpha=(k_1,k_2,\ldots,k_\ell)\in T_\infty$ we define
  \begin{equation*}
    \xt_\alpha=\begin{cases}  x_{(k_{n_1},k_{n_2},\ldots,k_{n_i})}
    &\text{if $\ell=n_i$ for some $i\in\N$}\\
    z_{k_\ell}        &\text{if $\ell\in \N\setminus N$}.
    \end{cases}
  \end{equation*}
  Our claim now follows from the fact that
  $(\xt_\alpha)_{\alpha\in T_\infty}$ is also a  normalized, weakly
  null tree and that for any branch $(y_i)$ of
  $(\xt_\alpha)_{\alpha\in T_\infty}$ the subsequence $(y_{n_i})$ is a
  branch of $(x_\alpha)_{\alpha\in T_\infty}$.
\end{proof}

In the definition of $U$-upper and $V$-lower tree estimates it is
actually not necessary to assume that $C$ exists uniformly for all trees
as the following proposition shows.

\begin{prop}\label{P:1.2a}
  Let $U$ and $V$ be Banach spaces with normalized, 1-unconditional
  bases $(u_i)$ and $(v_i)$. Assume that $X$ is a Banach space with
  the property that every normalized, weakly null tree in $X$ has a
  branch which dominates $(v_i)$ and a branch which is dominated by
  $(u_i)$. Then $X$ satisfies $(V,U)$-tree estimates.
\end{prop}
\begin{proof}
  Assume that $X$ has the property
  that for any $C\ge 1$ there is a weakly null tree in
  $S_X$ so that every branch  does not $C$-dominate $(v_i)$ (the
  argument for $U$ is similar). We will choose a tree
  $(x_\alpha)_{\alpha\in T_\infty}$
  which has the property that each branch of $(x_\alpha)_{\alpha\in
  T_\infty}$ does not dominate $(v_i)$.

  By induction we will choose for every $m\in\N$ a well-founded subset
  $S_m\subset T_\infty$
  and a family $(x_\alpha^{(m)})_{\alpha\in S_m}$ in $S_X$ so that
  \begin{enumerate}
  \item[a)]$S_{m-1}\subset S_m$ (with $S_0=\emptyset$), and
    $\max(S_{m-1})\cap\max(S_m)=\emptyset$;
  \item[b)] for any $\alpha=(n_1,\ldots,n_\ell)\in
  S_m\cup\{\emptyset\}$ either $\alpha$ is maximal in $S_m$ or  $S_m$
  contains all the direct successors $(\alpha,n)$, $n>n_\ell$ (put
  $\ell=n_\ell=0$ if $\alpha=\emptyset$), of $\alpha$;
  \item[c)] $(x^{(m)}_{(\alpha,n)})_{n>n_\ell}$ is a weakly null sequence
    for any $\alpha=(n_1,\ldots,n_\ell)\in S_m\setminus \max(S_m)$;
  \item[d)] if $\alpha=(n_1,\ldots,n_\ell)\in \max(S_m)$
    and $k\in\{1,2,\dots,\ell\}$ is such that
    $(n_1,n_2,\ldots,n_k)\in\max(S_{m-1})$ ($k=0$ if $m=1$), then the
    segment $(y_i)_{i=k+1}^\ell$, where
    $y_i=x^{(m)}_{(n_1,n_2,\ldots,n_i)}$ for $i=k+1,k+2,\ldots,\ell$,
    does not $m$-dominate $(v_i)_{i=k+1}^\ell$.
  \end{enumerate}
  Once we have finished the construction of $S_m$ and
  $(x^{(m)}_\alpha)_{\alpha\in S_m}$ we
  deduce from~(a) and~(b) that $\bigcup_{m\in\N} S_m=T_\infty$.
  For $\alpha\in T_\infty$ let $x_\alpha=x^{(m)}_\alpha$ with
  $m=m(\alpha)=\min\{m': \alpha\in S_{m'}\}$.
  Let $(y_n)$ be a branch of $(x_\alpha)_{\alpha\in T_\infty}$, say
  $y_n=x_{\alpha_n}$ for $n\in\N$ and for
  some increasing (with respect to extension) sequence
  $(\alpha_n)\subset T_\infty$ with $|\alpha_n|=n$ for $n\in\N$.
  For $m\in\N$ let $\ell_m=\max\{ \ell: \alpha_\ell\in S_m\}$ and deduce
  from (d) that
  $(y_i)_{i=\ell_{m-1}+1}^{\ell_m}= (x^{(m)}_{\alpha_i}
  )_{i=\ell_{m-1}+1}^{\ell_m}$
  does not  $m$-dominate $(v_i)_{i=\ell_{m-1}+1}^{\ell_m}$. Thus $(y_n)$
  does not dominate $(v_n)$.
  
  Assume we have chosen $S_{m-1}$ and $(x^{(m-1)}_\alpha)_{\alpha\in
  S_{m-1}}$ for some $m$.
  For $\alpha \in\max( S_{m-1})$ we can choose a  normalized  weakly null tree
  $(z^{(\alpha)}_\beta)_{\beta\in T_\infty}$ (if $m=1$, and thus
  $S_{m-1}=\emptyset$, we
  choose one tree    $(z_\beta)_{\beta\in T_\infty}$)
  so that no branch  $3(|\alpha|+m)$-dominates $(v_i)$. Since
  $(z^{(\alpha)}_\beta)_{\beta\in T_\infty}$  is weakly null we can,
  after passing to an
  appropriate
  full subtree,  assume that every branch of
  $(z^{(\alpha)}_\beta)_{\beta\in T_\infty}$ is a basic sequence with
  projection constant less than $3$.
  
  For $\beta=(n_1,n_2,\ldots,n_\ell)\in T_\infty\setminus S_{m-1}$
  we define $\alpha(\beta)$ to be the maximal restriction of $\beta$
  which lies in $S_{m-1}$
  and  let $y^{(\beta)}_i=z^{(\alpha(\beta))}_{(n_1,n_2,\ldots,n_i)}$ for
  $i=1,\ldots,\ell$. Define 
  \begin{equation*}
    S_m=  S_{m-1}\cup
    \left \{ \beta\in T_\infty\setminus S_{m-1}:
    (y_i^{(\beta)})_{i=1}^{|\beta|-1}\ \text{$3(|\alpha(\beta)|+m)$-dominates }
    (v_i)_{i=1}^{|\beta|-1}
    \right\}.
  \end{equation*}
  $S_m$ is well-founded, otherwise for some $\alpha\in\max(S_{m-1})$
  some branch of $(z^{(\alpha)}_\beta)_{\beta\in T_\infty}$ would
  $3(|\alpha|+m)$-dominate $(v_i)$. Since for any $\beta\in
  T_\infty\setminus S_{m-1}$ the sequence
  $(y^{(\beta)}_i)_{i=1}^{|\alpha(\beta)|}$
  $3|\alpha(\beta)|$-dominates $(v_i)_{i=1}^{|\alpha(\beta)|}$~(a)
  holds.

  If $\beta\in S_m\setminus S_{m-1}$ is not maximal in $S_m$, then 
  $(y^{(\beta)}_i)_{i=1}^{|\beta|}$ $3(|\alpha(\beta)|+m)$-dominates
  $(v_i)_{i=1}^{|\beta|}$. It follows that any direct successor of
  $\beta$ is in $S_m$, which implies condition (b). Finally we put for
  $\gamma\in S_m$
  $$x^{(m)}_\gamma=\begin{cases}  x^{(m-1)}_\gamma &\text{if
    $\gamma\in S_{m-1}$}\\
  \phantom{----}\\
  z^{(\alpha(\gamma))}_\gamma
  &\text{if $\gamma\not\in S_{m-1}$.}
  \end{cases}$$
  Condition (c) is satisfied and, for $\beta\in \max(S_m)\subset
  T_\infty\setminus S_{m-1}$,  if
  $(y_i^{(\beta)})_{i=|\alpha(\beta)|+1}^{|\beta|}$ 
  $m$-dominates  $(v_i)_{i=|\alpha(\beta)|+1}^{|\beta|}$,
  then $(y_i^{(\beta)})^{|\beta|}_{i=1}$
  $3(|\alpha(\beta)|+m)$-dominates $(v_i)_{i=1}^{|\beta|}$, which is
  not true. Thus (d) holds.
  
\end{proof}
The following Proposition generalizes a result of Prus \cite{P}.
\begin{prop}
  \label{P:1.3}
  Assume that $Z$ has an FDD $(E_i)$, and let $V$ be a space with
  a normalized and 1-unconditional basis $(v_i)$.
  
  The following statements are equivalent:
  \begin{enumerate}
  \item[a)] $(E_i)$ satisfies $V$-lower estimates in $Z$,
  \item[b)] $(E^*_i)$ satisfies $\Vs$-upper estimates in  $\Zs$.
  \end{enumerate}
  (Here $\Vs$-upper estimates are with respect to $(v_i^*)$, the
  sequence of biorthogonal functionals to $(v_i)$).
  
  Moreover, if $(E_i)$ is bimonotone in $Z$, then the equivalence
  holds true if one replaces, for some $C\ge 1$, $V$-lower estimates
  by $C$-$V$-lower estimates in (a) and  $\Vs$-upper estimates by
  $C$-$\Vs$-upper estimates in (b).
\end{prop}
\begin{remark} By duality Proposition~\ref{P:1.3} holds if we
  interchange the words \emph{lower} and \emph{upper} in~(a) and~(b).
\end{remark}

\begin{proof} W.l.o.g. we assume that $E=(E_i)$ is bimonotone in $Z$.

\noindent``(a)$\Rightarrow$(b)'' Assume that $(E_i)$ satisfies $C$-$V$-lower
 estimates in $Z$, and let  $(z_i^*)_{i=1}^\ell$ be a block sequence of
 $E^*=(E_i^*)$.

 For an appropriate $z\in S_Z$
 with $\supp_E(z)\subset
 [\min\supp_{E^*}(z_1^*),\max\supp_{E^*}(z_l^*)]$ we have (putting
 $\max\supp_{E^*}(z_{0}^*)=0$)
\begin{align*}
\Bigl\|\sum_{i=1}^\ell z_i^*\Bigr\|_{Z^{(*)}}
&=\sum_{i=1}^\ell z_i^*(z)\\
&= \sum_{i=1}^\ell z_i^*
   \big(P^E_{(\max\supp_{E^*}(z_{i-1}^*),\max\supp_{E^*}(z_{i}^*)]}(z)\big)\\
&\leq \sum_{i=1}^\ell \|z_i^*\|\cdot
 \|P^E_{(\max\supp_{E^*}(z_{i-1}^*),\max\supp_{E^*}(z_{i}^*)]}(z)\|\\
&\leq\Bigl\|\sum_{i=1}^\ell \|z_i^*\|  \cdot v_i^*\Bigr\|_{\Vs}
  \Bigl\|\sum_{i=1}^\ell  \|P^E_{(\max\supp_{E^*}(z_{i-1}^*), \max\supp_{E^*}(z_{i}^*)]}(z)
 \| v_i\Bigr\|_{V} \\
&\leq C\Bigl\|\sum_{i=1}^\ell \|z_i^*\|  \cdot v_i^*\Bigr\|_{\Vs}
  \Bigl\|\sum_{i=1}^\ell  \|P^E_{(\max\supp_{E^*}(z_{i-1}^*), \max\supp_{E^*}(z_{i}^*)]}(z)  \|
 z_i\Bigr\| \\
&\text{\Big(where }
 z_i=\frac{P^E_{(\max\supp_{E^*}(z_{i-1}^*), \max\supp_{E^*}(z_{i}^*)]}
         (z)}{\|P^E_{(\max\supp_{E^*}(z_{i-1}^*), \max\supp_{E^*}(z_{i}^*)]}(z)  \|}
        \text{ for $i=1,\ldots,\ell$} \text{ with }\frac00=0\Big)\\
&= C\Bigl\|\sum_{i=1}^\ell \|z_i^*\| \cdot v_i^*\Bigr\|_{\Vs}  \|z\|
= C\Bigl\|\sum_{i=1}^\ell \|z_i^*\| \cdot v_i^*\Bigr\|_{\Vs}.
\end{align*}
This implies (b).

\noindent``(b)$\Rightarrow$(a)'' Assume that $(E_i^*)$ satisfies $C$-$\Vs$-upper
 estimates in $\Zs$, and let  $(z_i)_{i=1}^\ell$ be a block sequence of
 $(E_i)$.
Choose $\sum_{i=1}^\ell a_i v_i^*\in S_{\Vs}$ so that
 $$\sum_{i=1}^\ell a_i\|z_i\|=
 \Bigl(\sum_{i=1}^\ell a_iv^*_i\Bigr)\Bigl(\sum_{i=1}^\ell \|z_i\| v_i\Bigr)=
                                                   \Bigl\|\sum_{i=1}^\ell \|z_i\|
                   v_i\Bigr\|_V,$$
and choose, for $i=1,2,\ldots,\ell$, a vector $z_i^*\!\in\! S_{\Zs}$ with
  $z_i^*(z_i)=\|z_i\|$ and
 $$\supp_{E^*} (z_i^*)\!\subset\![\min\supp_E (z_i), \max\supp_E (z_i)].$$
 ($(E_i)$ is assumed to be bimonotone in $Z$.) Let $z^*=\sum_{i=1}^\ell a_i z_i^*$. It
 follows from our assumption (b) that $\|z^*\|\le C$, and thus that
$$\Bigl\|\sum_{i=1}^\ell z_i\Bigr \|\ge    \frac1C \sum_{i=1}^\ell z^*(z_i)=
 \frac1C \sum_{i=1}^\ell a_i\|z_i\|=\frac1C \Bigl\|\sum_{i=1}^\ell \|z_i\| v_i\Bigr\|_V,$$
which implies (a).
\end{proof}
\begin{prop}\label{P:1.4} Assume that $U$ is a space
  with a normalized and 1-unconditional basis $(u_i)$ and that $X$ is a reflexive space which
  satisfies  $C$-$U$-upper tree estimates for some $C\ge 1$.

Then, for any $\vp>0$, $X^*$ satisfies $(2C+\vp)$-$\Us$-lower tree estimates.
\end{prop}
\begin{remark} One might ask, whether or not
the converse of Proposition \ref{P:1.4} is true, i.e., similar to the
FDD case, whether
 $X$ satisfies $U$-upper tree estimates if
 $X^*$ satisfies $\Us$-lower tree estimates.

The answer is affirmative under certain conditions on $U,$ but we do
not give a direct proof for that fact. Once we  have shown
that, under appropriate conditions, a separable, reflexive space
 $X$ which satisfies $U$-lower tree estimates is both a subspace and a
quotient of spaces
 having an
 FDD   satisfying $U$-lower estimates,  this result  will
 follow easily
 (see Corollary \ref{C:3.2a} in Section \ref{S:3}).
\end{remark}
\begin{proof} Let $\eta>0$
  and let $(x^*_i)$ be a normalized, weakly null sequence in $X^*$.
Then there is a subsequence  $(x^*_{i_n})$  and a normalized, weakly null
 sequence $(y_n)$ in $X$ so that $x^*_{i_n}(y_n)>\frac1{2+\eta}$ for $n\in\N$.
 Indeed, for each $i\in\N$ choose $x_i\in S_X$ with $x^*_i(x_i)=1$, take a subsequence
   $(x_{j_n})$ so that $x=w$-$\lim_{n\to\infty} x_{j_n}$ exists, and
then for an appropriately large $n_0$ and for each $n\in\N$ let
$$y_n= \frac{x_{j_{n+n_0}}-x}{\|x_{j_{n+n_0}}-x\|}$$
and $i_n=j_{n+n_0}$.

Now let $(x^*_\alpha)_{\alpha\in T_\infty}$ be a normalized, weakly null
 tree in $X^*$. By replacing  certain nodes by subsequences,
 using the previous observation,
 we can pass to a full subtree
   $(\tilde x^*_\alpha)_{\alpha\in T_\infty}$    for which
there is a normalized, weakly null tree $(y_\alpha)_{\alpha\in T_\infty}$ with
 $\tilde x^*_\alpha(y_\alpha)>\frac1{2+\eta}$ for all $\alpha\in T_\infty$.
 Secondly, we may assume, again after passing to full  subtrees, that
$|\tilde x^*_\alpha(y_\beta)|<2^{-m-n}\eta$
and  $|\tilde x^*_\beta(y_\alpha)|<2^{-m-n}\eta$
 whenever $\alpha,\beta\in T_\infty$, $|\alpha|=m<|\beta|=n$
  and $\beta$ is an extension of $\alpha$.

By our assumption we can extract  a branch  $ (z_n)$ from
$(y_\alpha)_{\alpha\in T_\infty}$ which is $C$-dominated by $(u_i)$.
Let $ (z_n^*)$ be the corresponding branch of
$(\tilde x^*_\alpha)_{\alpha\in T_\infty}$, and let
$(a_i)\in\coo$. Choose $(b_i)\in \coo$
 with $\|\sum b_i u_i\|=1$ and
$\sum a_i b_i =(\sum a_i u_i^*)(\sum b_i u_i)=\|\sum a_i u_i^*\|$.
It follows that $\|\sum b_i z_i\|\le C$, and thus that (note that $a_ib_i\ge 0$ for $i\in\N$)
\begin{align*}
\Bigl\| \sum a_i z_i^*\Bigr\|&\ge \frac1C
\sum a_i z_i^*\Big( \sum b_i z_i    \Big)\\
&\ge  \frac1C \sum a_i b_iz^*_i(z_i)-\frac1C\sum_{i\not= j}
|a_ib_j|\cdot|z^*_i(z_j)|
 \ge
\frac1C \frac1{2+\eta}\sum a_ib_i - \frac1C\max_{i,j} |a_ib_j|\eta,
\end{align*}
which implies our claim if we choose $\eta>0$ small enough.
\end{proof}

The following connection between lower and upper tree
estimates and lower and upper estimates for spaces with FDDs
will be shown with techniques  developed in  \cite{OS1} and \cite{KOS}.
\begin{prop}\label{P:1.5}
  Assume that $V$ is a Banach space with a normalized  and
  1-unconditional basis $(v_i)$, and let $Z$ be a reflexive space with
  an FDD $(E_i)$.

  If $Z$ satisfies $V$-upper or -lower tree estimates, then $(E_i)$
  can be blocked into an FDD $(F_i)$ which satisfies   $V$-upper or
  -lower   estimates in $Z$.
\end{prop}

For the proof of Proposition \ref{P:1.5} we will need to recall some
notation and a proposition from~\cite{OS1} and~\cite{OS2}.

\begin{defin}
If $\cA \subseteq S_X^\omega$, the set of all normalized sequences in $X$, and
$\vp >0$, we set
$$\cA_\vp = \{(x_n)\in S_X^\omega :\text{ there exists } (y_n) \in \cA
\text{ with } \|x_n-y_n\| < \frac{\vp}{2^n}\text{ for all } n\}\ .$$
\end{defin}

\noindent
$\overline{\cA_\vp}$ denotes the closure
of $\cA_\vp$ w.r.t.\ the product topology of the
discrete topology on $S_X$.

The next result follows from Proposition~2.4 in~\cite{OS2} (which is a
restatement of a part of Theorem~3.3 in~\cite{OS1}) and
Proposition~2.5 in \cite{OS2}. In this section we will only use it for
the (much simpler) case $X=Z$. In section \ref{S:3} we will use it in
its full generality.
\begin{prop}\label{P:3.3}
  Let $X$ be a Banach space which is a subspace of a reflexive space
  $Z$ with an FDD $(E_i)$. Let $\cA\subseteq S_X^\omega$. Then the
  following are equivalent.
  \begin{enumerate}
  \item[{a)}] For all $\vp>0$ every normalized, weakly null tree in
    $X$ has a branch in $\overline{\cA_\vp}$.
  \item[{b)}]  For all $\vp>0$ there exists a blocking $(F_i)$ of $(E_i)$
    and $\deltab = (\delta_i)$, $\delta_i\downarrow 0$, so that if
    $(x_n) \subseteq S_X$ is a $\deltab$-skipped block sequence of
    $(F_i)$ in $Z$,
    then $(x_n) \in \overline{\cA_\vp}$.
  \end{enumerate}
\end{prop}

\begin{proof}[Proof of Proposition \ref{P:1.5}]
If $Z$ satisfies $V$-upper tree estimates, then $Z^*$ satisfies, by
Proposition \ref{P:1.4},  $\Vs$-lower tree estimates, and if we can
block $(E^*_i)$ into an FDD $(H^*_i)$ which
 satisfies $\Vs$-lower estimates in $Z^*$,  then,  by Proposition
 \ref{P:1.3} and the remark following it,
 $(H_i)$ satisfies $V$-upper estimates in $Z$.

 Therefore we need to prove the Proposition only for the case that $Z$
 satisfies $V$-lower  tree estimates.

 Let $K$ be the projection constant of $(E_i)$ in $Z$, and choose
 $C\geq 1$ such that the space $Z$ satisfies  $C$-$V$-lower tree
 estimates. Applying the implication ``(a)$\Rightarrow$(b)'' of
 Proposition~\ref{P:3.3} to $X=Z$,
 $$\cA=\Big\{ (z_i)\in S_Z^\omega:\,(z_i)\ C\text{-dominates
 }(v_i)\Big\}$$
 and to an $\vp>0$ small enough so that
 $$\overline{\cA_\vp}\subset\Big\{ (z_i)\in S_Z^\omega:\, (z_i)\
 2C\text{-dominates }(v_i)\Big\}\ ,$$
 we obtain a blocking $(F_i)$ of $(E_i)$  and a sequence
 $\deltab'=(\delta_i')$, $\delta_i'\downarrow 0$, so that
 every $\deltab'$-skipped block sequence of $(F_i)$ in $Z$
 $2C$-dominates $(v_i)$. By the remark following the definition of a
 $\deltab$-skipped block sequence we may assume, after replacing
 $\deltab'$ by $\frac1{2K}\deltab'$ if necessary, that in fact every
 $\deltab'$-skipped block sequence of any subsequent blocking of
 $(F_i)$ in $Z$ $2C$-dominates $(v_i)$.

 By Proposition~\ref{tree-est} $Z$ also
 satisfies $C$-$[v_{i+1}]_{i=1}^{\infty}$-lower tree estimates. Hence we
 can repeat the above argument to obtain a further blocking $G=(G_i)$
 of $(F_i)$ and a sequence $\deltab=(\delta_i)$, $\delta_i\downarrow 0$
 and $\delta_i\leq
 \delta_i'$ for all $i$, so that every $\deltab$-skipped block
 sequence of $(G_i)$ in $Z$ $2C$-dominates $(v_i)$
 \emph{and}
 $(v_{i+1})$. W.l.o.g.~we can assume that
 $\sum_{i=1}^\infty\delta_i\le \frac1{32C}$.
 Using a result in~\cite{J} (see also~\cite{KOS}, Lemma~4.2) we can block
 $(G_i)$ into $H=(H_i)$, say $H_i=\oplus_{j=N_{i-1}+1}^{N_i} G_j$ for
 $i\in\N$, where $0=N_0<N_1<N_2<\ldots$, so that
 for any $z\in S_Z$, $z=\sum_{j=1}^\infty x_j$ with $x_j\in G_j$ for $j\in\N$,
 and for every $i\in\N$ there is a $t_i\in(N_{i-1},N_i)$ so that
 $$\|P^G_{t_i}(z)\|=\|x_{t_i}\|<\delta_i^2.$$
 Assume now that $(z_n)$ is a normalized block sequence of $(H_i)$ in $Z$.
 We will show that
 $$\Big\|\sum a_i z_i\Big\|\ge \frac1{16C}\qquad
 \text{whenever }\Big\|\sum a_i v_i\Big\|=1.$$
 
 For $i\in\N$ let $k_i\in\N$ such that
$z_i\in \oplus_{j=k_{i-1}+1}^{k_i} H_j$ ($k_0=0$) and choose some
$t_i \in (N_{k_i-1},N_{k_i})$
 for which $\|P^G_{t_i}(z_i)\|<\delta_i^2$.  For $i\in\N$ write
 $z_i=g_i+h_i$, where
$$g_i=P^G_{[N_{k_{i-1}}+1,t_i]}(z_i)\text{ and } h_i=P^G_{(t_i,N_{k_{i}}]}(z_i) ,$$
and let
$$B=\big\{i\in\N: \|a_ih_i +a_{i+1} g_{i+1}\|\ge \delta_i \big\}.$$
 We define for $i\in\N$
$$w_i = \frac{a_ih_i +a_{i+1} g_{i+1}}{\|a_ih_i +a_{i+1} g_{i+1}\|}$$
 ($0/0=0$) and
\begin{align*}
&\alpha_i=\begin{cases} \|a_ih_i +a_{i+1} g_{i+1}\| &\text{if $i\in B$}\\
                                       0 &\text{if $i\not\in B$.}
                \end{cases}
\end{align*}
For $i\in\N$ we put $\wt_i=w_i$ if $i\in B$ and we let $\wt_i$
 be some normalized element in $G_{N_{k_i}}$ if $i\in \N\setminus B$.
Note that $(\tilde w_i)$ is a $\deltab$-skipped block sequence of $(G_i)$, and
we deduce that
\begin{align*}
  \Big\|\sum &a_i z_i\Big\|=\Big\|\sum a_ig_i+a_ih_i\Big\|\\
  &=\Big\| a_1 g_1+\sum_{i=1}^\infty  \|a_ih_i +a_{i+1} g_{i+1}
  \|w_i\Big\|\\
  &\ge\Big\| a_1 g_1+\sum_{i\in B}  \|a_ih_i +a_{i+1} g_{i+1}\|
  w_i\Big\|-\frac1{32C}\\
  &\ge \frac12\Big[|a_1| \|g_1\|+ \Big\|\sum_{i\in B}  \|a_ih_i
    +a_{i+1} g_{i+1}\|w_i\Big\|\Big]-\frac1{32C}\\
  &= \frac12\Big[|a_1| \|g_1\|+ \Big\|\sum_{i\in\N}      \alpha_i
    \tilde w_i\Big\|    \Big]-\frac1{32C}\\
  &\ge  \frac1{8C}\Big[|a_1| \|g_1\|+ \Big\|\sum_{i\in\N}
    \alpha_i v_i \Big\|+ \Big\|\sum_{i\in\N}
    \alpha_i v_{i+1} \Big\|   \Big]-\frac1{32C}\\
  &\geq\frac1{8C}\Big[|a_1| \|g_1\|+ \Big\|\sum_{i\in\N}   \|a_ih_i
    +a_{i+1} g_{i+1}\|   v_i \Big\|+ \Big\|\sum_{i\in\N}   \|a_ih_i
    +a_{i+1} g_{i+1}\|   v_{i+1} \Big\|   \Big]-\frac1{16C}\\
  &\geq\frac1{8C}\Big[|a_1| \|g_1\|+ \Big\|\sum_{i\in\N}   \|a_ih_i\|
  v_i \Big\|+ \Big\|\sum_{i\in\N}   \|a_{i+1} g_{i+1}\|   v_{i+1}
  \Big\|   \Big]-\frac1{16C}\\
  &\ge  \frac{1}{8C}\Big\|\sum_{i\in\N}   \|a_ih_i+a_ig_i\|v_i\Big\|
  -\frac1{16C}\\
  &= \frac1{8C}\Big\|\sum_{i\in\N}   a_iv_i\Big\| - \frac1{16C}  =
  \frac1{16C},
\end{align*}
which finishes the proof of our claim.\end{proof}

\section{The space $Z_V(E)$}\label{S:2}

Let $Z$ be a space with an FDD $E=(E_i)$, and let $V$ be a space with a
1-unconditional and normalized
basis $(v_i)$. The space $Z_V=Z_V(E)$ is defined to be the completion
of $c_{00}(\oplus E_i)$ with respect to the following norm $\|\cdot\|_{Z_V}$.
$$\|z\|_{Z_V}
=\max_{\stackrel{k\in\N}{0=n_0<n_1<n_2<\ldots<n_k}}
\Big\|\sum_{j=1}^k\|P^E_{(n_{j-1},n_j]}(z)\|_Z\cdot
  v_j\Big\|_V\qquad\text{for }z\in  c_{00}(\oplus E_i).$$
Note that $(E_i)$ is a monotone FDD in $Z_V$, which implies that
the projection constant of $(E_i)$ in $Z_V$ is at most $2$.
$(E_i)$ is bimonotone in $Z_V$ if it is bimonotone  in $Z$, and if
$Z$ and $Z'$ are isomorphic
and $E'=(E'_i)$ is the image of $E$ under
an isomorphism,
then $Z'_V(E')$ and $Z_V(E)$ are naturally isomorphic.

\begin{lem}\label{L:2.1} Assume that $V$ is a Banach space with   a normalized
 and $1$-unconditional basis $(v_{i})$ and  that
 $(v_{i+1})$  $C$-dominates $(v_i)$ for some $C\ge1$. Let $Z$ be a space with an FDD $E=(E_i)$.

Then every normalized block sequence $(z_n)$ of $(E_i)$ in
 $Z_V(E)$ has a subsequence $(z_{i_n})$  for which there is a normalized
 block sequence $(b_n)$ in $V$ so that for some $d>0$
\begin{equation}\label{E:2.2.1}
\Bigl\|\sum_{i=1}^\infty a_iz_i\Bigr\|_{Z_V}\ge d
\Bigl\|\sum_{n=1}^\infty a_{i_n} b_{n}\Bigr\|_V\text{ whenever }(a_i)\in \coo.
 \end{equation}
In particular $(z_{i_n})$ dominates $(b_n)$ (choose $a_i=0$ if $i\not\in\{i_1,i_2,\ldots\}$).
\end{lem}
\begin{proof} We assume without loss of generality that $(E_i)$  is
 bimonotone in $Z$.
 Let $(z_i)$ be  a normalized block sequence of $(E_i)$ in $Z_V$.
We will choose a subsequence
$(z_{i_j})_{j=1}^\infty$ together
with $\vp_0>0$ and increasing sequences $(m_j)_{j=1}^\infty$ and
$(n_s)_{s=1}^\infty$ in $\N$ so that for all $j\in \N$
\begin{eqnarray}
  \label{E:2.2.2}
  & n_{m_j}=\max\supp_E(z_{i_j}), \text{ and }&\\[3ex]
  \label{E:2.2.3}
  & \Bigl\|\sum_{s=m_{j-1}+1}^{m_j}\|
  P^E_{(n_{s-1},n_s]}(z_{i_j})\|_Z\cdot v_s\Bigr\|_V\ge \vp_0 &
    \text{\hspace{1em}(where $m_0=n_0=0$).}
\end{eqnarray}
Then (\ref{E:2.2.1}) follows immediately with
$b_j=\tilde b_j/\| \tilde b_j\|$ and
$$\tilde b_j=\sum_{s=m_{j-1}+1}^{m_j}\|
P^E_{(n_{s-1},n_s]}(z_{i_j})\|_Z\cdot v_s\qquad\text{for $j\in \N$.}$$
Indeed, if $(a_i)\in c_{00}$ and $z=\sum a_i z_i$, then
\begin{align*}
\|z\|_{Z_V}&\ge \Big\| \sum_{j=1}^\infty \sum_{s=m_{j-1}+1}^{m_j}
\|P^E_{(n_{s-1},n_s]}(z)\|_Z\cdot v_s\Big\|_V\\
  & \ge \Big\| \sum_{j=1}^\infty \sum_{s=m_{j-1}+1}^{m_j} |a_{i_j}|\cdot
  \|P^E_{(n_{s-1},n_s]}(z_{i_j})\|_Z\cdot v_s\Big\|_V\\
&\quad\text{(using bimonotonicity)}\\
&= \Big\| \sum a_{i_j} \tilde b_j\Big\|_V
\ge  \vp_0\Big\| \sum a_{i_j} b_j  \Big\|_V.
\end{align*}

We can, for each $i\in\N$, choose positive integers
$k(i)$ and $0=n_0{(i)}<n_1{(i)}<\ldots<n_{k(i)}{(i)}=\max\supp_E(z_i)$ so that
\begin{equation}\label{E:2.2.3a}
  \|z_i\|_{Z_V}=
  \Bigl\|\sum_{j=1}^{k(i)}\| P^E_{(n_{j-1}{(i)},n_j{(i)}]}(z_i)\|_Z
    \cdot v_j \Big\|_V=1.
\end{equation}
We can assume that we are in one of the following three cases:

\noindent Case 1: $\|z_{i_n}\|_Z\ge \vp_0$ for all $n\in\N$, some
$\vp_0>0$ and some subsequence $(z_{i_n})$ of $(z_i)$.

\noindent Case 2: $\bigl\|\sum_{j=1}^{\max\supp_E(z_{i_n})} \|
P^E_{j}(z_{i_n})\|_Z \cdot v_j\bigr\|_V\ge \vp_0$ for all $n\in\N$,
some $\vp_0>0$ and some subsequence $(z_{i_n})$ of $(z_i)$.

\noindent Case 3:  $\lim_{i\to\infty}\|z_i\|_Z=0$ and
$$\lim_{i\to\infty} \Bigl\|\sum_{j=1}^{\max\supp_E(z_i)}
\| P^E_{j}(z_i)\|_Z\cdot v_j\Bigr\|_V=0.$$

\noindent
Indeed, if all subsequences $(z_{i_n})$ fail Cases~1 and~2, then Case~3
holds.

In Case~1 we choose $n_j=\max\supp(z_{i_j}),\ m_j=j$, and in case~2 we
choose $n_j=j,\ m_j=\max\supp(z_{i_j})$ for each $j\in\N$. In case 3
we will choose by induction $i_j$, $m_j$ and $n_{m_{j-1}+1}<
n_{m_{j-1}+2}< \ldots <n_{m_j}$, $j\in\N$, so that $(i_j)$, $(m_j)$
and $(n_s)$ are increasing and so that~(\ref{E:2.2.2})
and~(\ref{E:2.2.3})  are satisfied with $\vp_0=1/2$  for  all
$j\in\N$.

For $j=1$ we choose $i_1=1$, $m_1=k(1)$ and $n_s=n_s(1)$ for
$s=1,2,\ldots, m_1$.
Assume we have chosen $i_1<i_2<\ldots<i_{j-1},\
m_1<m_2<\ldots<m_{j-1}$ and $n_1<\ldots< n_{ m_{j-1}}$.  By the
first condition of Case 3 we can choose $i'>i_{j-1}$ large enough so
that for all  $i\ge i'$ there is a $k'(i)\in (m_{j-1},k(i))$ with (we
are using that $\max_{s\le k(i)} \|P^E_{(n_{s-1}(i),n_s(i)]}
  (z_i)\|_Z\to 0$ as $i\to\infty$)
\begin{equation}\label{E:2.2.4}
1/3<
\Bigl\|\sum_{s=1}^{k'(i)}\|
P^E_{(n_{s-1}{(i)},n_s{(i)}]}(z_i)\|_Z\cdot v_{s}\Big\|_V<1/2,
\end{equation}
and thus, by (\ref{E:2.2.3a}),
\begin{equation}\label{E:2.2.5}
\Bigl\|\sum_{s=k'(i)+1}^{k(i)}
 \| P^E_{(n_{s-1}{(i)},n_s{(i)}]}(z_i)\|_Z \cdot v_{s}\Big\|_V>1/2.
\end{equation}
Since by our assumption $(v_{s_0+s})_{s=1}^\infty$ $C^{s_0}$-dominates
$(v_s)_{s=1}^\infty$ for all $s_0\in\N$, it follows, using the second
condition in Case 3, that $n_{k'(i)}(i)-k'(i)\to\infty$ as $i\to\infty$.
Indeed, assuming this is not true, we can choose an infinite subset
$N\subset\N$ so that
for some $M\ge 0$ we have $n_{k'(i)}(i)-k'(i)=M$ for all $i\in N$. It
follows that for each $i\in N$ at most $M$ of the intervals
$(n_{s-1}(i),n_s(i)]$ in~\eqref{E:2.2.4} has length
exceeding~$1$ and that the sum of the lengths of these intervals is at
most~$2M$, and hence their contribution to the norm
in~\eqref{E:2.2.4} converges to~$0$ as $i\to\infty$. Thus by a further
stabilization, replacing $N$ by an infinite subset of $N$ if
necessary, we may assume that for some $s_0\le M$ and for all $i\in N$
there is an interval $J_i\subset\{1,2,\ldots,k'(i)\}$ such that
$$\Big\|\sum_{j\in J_i} \|P^E_{j+s_0}(z_i)\|_Z\cdot
v_j\Big\|_V\ge \frac1{3(M+1)+1},$$
and thus, by our assumption on $(v_j)$, for some $\delta_0>0$
$$
\Big\|\sum_{j=1}^{\max \supp (z_i)} \|P^E_j(z_i)\|_Z\cdot
v_{j}\Big\|\geq 
\Big\|\sum_{j\in J_i} \|P^E_{j+s_0}(z_i)\|_Z\cdot v_{j+s_0}\Big\|\ge
\delta_0,$$
which contradicts the second assumption of Case 3.

Therefore we can choose $i_j=i>i'$ large enough so that
$n_{k'(i)}(i)-k'(i)>n_{m_{j-1}}$. Then set $m_j=k(i)$ and $n_s=n_s(i)$
for $k'(i)\leq s\leq k(i)$, choose
$n_{m_{j-1}+1}<n_{m_{j-1}+2}<\ldots< n_{k'(i)-1}$ arbitrarily from the
set $\big(n_{m_{j-1}},n_{k'(i)}(i)\big)$, and deduce our claim
from~(\ref{E:2.2.5}).
\end{proof}
\begin{cor}
  \label{C:2.1a}
  Assume that $(v_i)$ is a normalized, boundedly complete and
  1-unconditional basis of a Banach space $V$ so that $(v_{i+1})$
  dominates $(v_i)$, and let $Z$ be a space with an FDD $E=(E_i)$.

  Then $(E_i)$  is a boundedly complete FDD for $Z_V(E)$.
\end{cor}
\begin{proof}
  Let $(z_n)$ be a normalized block sequence of $(E_i)$ in
  $Z_V$. Choose a subsequence $(z_{i_n})$ of $(z_n)$ and a normalized
  block sequence $(b_n)$ in $V$ so that~\eqref{E:2.2.1} of
  Lemma~\ref{L:2.1} is
  satisfied for some $d>0$. Then it follows from the assumption that
  $(v_i)$ is boundedly complete that if $(a_i)\subset [\vp,\infty)$
  for some $\vp>0$, then
  $$\Bigl\|\sum_{j=1}^{i_n} a_j z_j\Bigr\|_{Z_V}\ge d
  \Bigl\|\sum_{j=1}^{n}a_{i_j} b_j\Bigr\|_V\to\infty\qquad\text{as
  }n\to\infty.$$
  Since $(z_i)$ was an arbitrary normalized
  block sequence of $(E_i)$ in $Z_V$
  it follows that $(E_i)$ is boundedly
  complete in $Z_V$.
\end{proof}

\begin{lem}\label{L:2.4}
Let $V$ be a Banach space with a normalized and $1$-unconditional basis $(v_i)$
 and assume that the space $Z$ has an FDD $E=(E_i)$.

    If $(v_i)$ is shrinking and if $(E_i)$ is shrinking in $Z$ then
 $(E_i)$ is shrinking in $Z_V(E)$ .
\end{lem}
\begin{proof} 
  W.l.o.g.~we assume that $(E_i)$ is bimonotone in $Z$ and therefore
  also in $Z_V(E)$. We first note that if $v^*=\sum_{i=1}^\infty a_i
  v^*_i$ converges in $V^*$ and  $\|v^*\|\le 1$ and if $(z_i^*)$ is a
  normalized block sequence of $(E^*_i)$ in $Z^*$, then the series
  $\sum_{i=1}^\infty a_i z^*_i$ converges in $(Z_V)^*$ and
  $\|\sum_{i=1}^\infty a_i z^*_i\|_{(Z_V)^*}\le 1$. Indeed, for $m\le n$
  in $\N$ there is a $z\in S_{Z_V}$ with
  $\supp_E(z)\subset[\min\supp_{E^*}(z_m^*),\max\supp_{E^*}(z_n^*)]$ so
  that
  \begin{align*}
    \Bigl\|\sum_{i=m}^n a_i z_i^*\Bigr\|_{(Z_V)^*}&=\sum_{i=m}^n a_i
    z^*_i(z)\\
    &\leq\sum_{i=m}^n |a_i|\cdot \|P^E_{(\max\supp_{E^*}
      (z^*_{i-1}),\max\supp_{E^*}(z^*_{i})]}(z)\|_Z\\
      &\le \Bigl\|\sum_{i=m}^n a_i v_i^*\Bigr\|_{V^*} \cdot
      \Bigl\|\sum_{i=m}^n\|P^E_{(\max\supp_{E^*}(z^*_{i-1}),
	\max\supp_{E^*}(z^*_{i})]}(z) \|_Z\cdot v_i \Bigr\|_{V}\\
	&\le\Bigl\|\sum_{i=m}^n a_iv_i^*\Bigr\|_{V^*}\cdot\|z\|_{Z_V}=
        \Bigl\|\sum_{i=m}^n  a_iv_i^*\Bigr\|_{V^*},
  \end{align*}
  which implies the claim.
  
  Define
  \begin{align*}
    K= \Bigl\{&\sum_{i=1}^\infty a_iz^*_i: \Bigl\|\sum_{i=1}^\infty
    a_i v^*_i\Bigr\|_{V^*}\le 1,\, (z_i^*)\text{ is an infinite block
      sequence in }S_{Z^*}\Bigr\}\\
    &\cup  \Bigl\{\sum_{i=1}^\ell a_i z^*_i:\ell\in\N,\,
    \Bigl\|\sum_{i=1}^\ell a_i v^*_i\Bigr\|_{V^*}\le 1,\,
    (z_i^*)\text{ is a block sequence in }S_{Z^*} \text{ of length
    }\ell\Bigr\},
  \end{align*}
  where  we allow in the second of the two sets which form $K$ the
  last element  $z^*_{\ell}$ to have infinite support.  Clearly, $K$
  is a $Z_V$-norming subset (isometrically) of $B_{(Z_V)^*}$. We claim
  that $K$ is $w^*$-compact. Indeed, let $y^*_n=\sum_{i=1}^\infty
  a^{(n)}_i z^*_{(n,i)}\in K$ for $n\in\N$ (where, for $n\in\N$,
  $a^{(n)}_i$ and $z^*_{(n,i)}$ may eventually vanish in case that
  $y^*_n$ is  in the second of the two sets which form $K$). After
  passing to a subsequence we can assume that
  \begin{align*}
    z^*_i&=\text{$w^*$-$\lim_{n\to\infty}$} z^*_{(n,i)} \in B_{Z^*}
    \text{ exists for all }i\in\N  \text{ and}\\
    v^*&=\text{$w^*$-$\lim_{n\to\infty}$} \sum a_i^{(n)} v^*_i  \in
    B_{V^*} \text{ exists.}
  \end{align*}
  Since $(v_i)$ is a shrinking basis of $V$, we can write $v^*=\sum
  a_i v_i^*$ for some scalars $(a_i)$. 
  
  Note that if $P^{E^*}_j(z_i^*)\not=0$, then
  $P^{E^*}_{j'}(z_{i'}^*)=0$ whenever $j'\le j$ and $i'>i$ or $j'\ge
  j$ and $i'<i$. This means that the non-zero terms of the sequence
  $(z^*_i)$ form a (finite or infinite) block sequence of
  $(E^*_i)$ (where, in the finite case, the last term may have
  infinite support).

  Since $\|z^*_i\|_{Z^*}\le 1$ for $i\in\N$, and since $(v_i)$ is
  1-unconditional, it follows that $$z^*=\sum_{i,\,
  z^*_i\not=0} a_i \|z_i^*\|_{Z^*}\frac{z_i^*}{\|z_i^*\|_{Z^*}}\in K.$$

  Finally, for $j\in\N$ and $z\in E_j$ we have
  $$\lim_{n\to\infty} y^*_n(z)=  \sum_{i=1}^\infty a_i
  z_i^*(z)=z^*(z),$$
  and thus $z^*$ is the $w^*$-limit of $(y_n^*)$. This shows that
  $K$ is $w^*$-closed.
  
  We deduce that $Z_V$ is embedded in $C(K)$, the space of continuous
  functions on $K$. Let $(z_i)$ be a bounded block sequence of $(E_i)$
  in $Z_V$, and let $z^*=\sum_{i=1}^\infty a_i z_i^* \in K$
  (i.e. $\|\sum_{i=1}^\infty a_i v_i^*\|_{V^*}\le 1$ and  either
  $\|z_i^*\|_{Z^*}=1$ for all $i\in\N$ or, for some $\ell\in\N$,
  $\|z_i^*\|_{Z^*}=1$ for all $i\le\ell$ and $z_i^*=0$ for $i>\ell$).
  If $(z^*_j)$ is an infinite normalized block sequence, then it
  follows that
  \begin{align*}
    z^*(z_i)&=\sum_{\stackrel{j\in\N,}{\max\supp_{E^*}(z_j^*)\ge
	\min\supp_E(z_i)}} a_j z^*_j(z_i)\\
    &\le \Bigl\|\sum_{\stackrel{j\in\N,}{\max\supp_{E^*}(z_j^*)\ge
	\min\supp_E(z_i)}} a_j v^*_j\Bigr\|_{V^*}\to0\qquad \text{as
	}i\to
    \infty.
  \end{align*}
  If, for some $\ell\in\N$, $\|z^*_\ell\|=1$ and $z^*_j=0$ for
  $j>\ell$, then from the assumption that $(E_i)$ is shrinking in $Z$
  and that $(z_i)$ is a bounded block sequence of $(E_i)$ in $Z_V$,
  and thus also in $Z$, we deduce that for large enough $i\in\N$
  $$z^*(z_i)=a_\ell z_\ell^*(z_i)\to 0\qquad\text{as }i\to \infty.$$
  It follows that $(z_i)$ is weakly null in $C(K)$, and thus in
  $Z_V$. Since $(z_i)$ was an arbitrary bounded block sequence in
  $Z_V$, this finishes the proof that $(E_i)$ is shrinking in $Z_V$.
\end{proof}
%\begin{remark}
 % The argument at the beginning of the proof shows the following. Let
 % $V$ be a Banach space with a normalized,
 % $1$-unconditional basis $(v_i)$, and let $Z$ be a Banach space with
 % an FDD $E=(E_i)$. Write $E^*$ for the FDD $(E_i^*)$ of
 % $Z^{(*)}$. Then $\big(Z_V(E)\big)^{(*)}$ is isometrically isomorphic to
 % $\big(Z^{(*)}\big)_{V^{(*)}}(E^*)$. We shall not need that in the
 % sequel.\marginpar{\scriptsize \hspace{4em}probably false!}
%\end{remark}
Corollary~\ref{C:2.1a} and Lemma~\ref{L:2.4} yield the following result.
\begin{cor}
  \label{C:2.4a}
  Assume that $Z$ is a space with a shrinking FDD $E=(E_i)$ and that
  $V$ is a reflexive Banach space with a normalized and
  $1$-unconditional basis $(v_i)$ such that $(v_{i+1})$ dominates
  $(v_i)$. Then $Z_V(E)$ is reflexive.
\end{cor}

We will now formulate  conditions on $V$ which ensure that,  given a
space $Z$ with an FDD $(E_i)$, every normalized block tree in
$Z_V(E)$ admits a branch that dominates some normalized block sequence
of $(v_i)$. We consider the following two  forms of shift invariance
of $V$.

\begin{defin}\label{D:2.2}
We say that $V$ has the {\em strong right  shift property} if
\begin{enumerate}
\item[(SRS)]  there exists $c>0$
 so that for all $(a_i)\in\coo$ and all $n\in\N$
$$\Bigl\|\sum_{i=1}^\infty a_i v_{i+n}\Bigr\|_V
\ge c\Bigl\|\sum_{i=1}^\infty a_i v_{i}\Bigr\|_V.$$
\end{enumerate}
We say that $V$ has the {\em weak left shift property} if
\begin{enumerate}
\item[(WLS)] there exists $d>0$ so that for all $m\in\N$ there exists
  $L=L(m)\ge m$ so that for all $k\le m$
$$\Bigl\|\sum_{i=L+1}^\infty a_i v_{i-k}\Bigr\|_V
\ge d\Bigl\|\sum_{i=L+1}^\infty a_i v_{i}\Bigr\|_V
 \text{ whenever }(a_i)\in\coo.$$
\end{enumerate}
\end{defin}
\begin{lem}\label{L:2.3}
 Let $Z$ be a space with an FDD $E=(E_i)$ and let $V$ be a space with
 a 1-unconditional and normalized basis $(v_i)$.

 If $V$ satisfies  $(SRS)$ and $(WLS)$, then there is a $C\ge 1$ so that
 every normalized block tree of $(E_i)$ in $Z_V(E)$  has a branch that
 $C$-dominates some normalized block sequence of $(v_i)$.
\end{lem}
\begin{proof} Without loss of generality we can, after renorming $Z$
 if necessary, assume that $(E_i)$ is bimonotone. Let $c$  and $d$ be
 as in  (SRS) and (WLS).

Given a block tree in $S_{Z_V(E)}$, we can extract a branch $(z_i)$ so that
\begin{multline}\label{E:2.3.1}
  L(b_{i-1})<a_i\quad\text{for all $i>1$, where}\ a_i=\min\supp(z_i)\
  \text{and}\\
  b_i=\max\supp(z_i)\ \text{for all $i\in\N$}.
\end{multline}
Using (SRS) and the fact that $(z_i)$ is normalized in $Z_V$,
we  can choose, for each $i\in\N$,
 $k(i)\in \N$ and $0=n_0{(i)}<n_1{(i)}<n_2{(i)}<\ldots<
n_{k(i)}{(i)}=b_i$ so that ($a_i\leq k(i)$ and)
\begin{align}\label{E:2.3.2}
  &n_j(i)=j\qquad\text{for } j=0,1,2,\ldots, a_i-1,\\
  &\label{E:2.3.3}  1=\|z_i\|_{Z_V}\ge
  \Bigl\|\sum_{j=1}^{k(i)}\| P^E_{(n_{j-1}{(i)},n_j{(i)}]}(z_i)\|_Z\cdot
  v_{j}\Bigr\|_V \ge c.
\end{align}
(Note that by  forcing (\ref{E:2.3.2}) we can only achieve the value
$c$.)

Put $m_1=k(1)$ and $n_j=n_j(1)$ for $j=0,1,2,\ldots, m_1$, and assume
that $m_1<m_2<\ldots<m_{i-1}$ and $n_0<n_1<\ldots<n_{m_{i-1}}=b_{i-1}$
have been chosen for some $i>1$. We put $m_i=m_{i-1}+k(i)-b_{i-1}$ and
$n_j=n_{j-m_{i-1}+b_{i-1}}(i)$ for $j=m_{i-1}+1, m_{i-1}+2,\ldots, m_i$.
Note that
\begin{gather*}
  m_i\geq m_{i-1}+a_i-b_{i-1}>m_{i-1},\\
  \text{\phantom{and}\qquad}
  n_{m_{i-1}}=b_{i-1}=n_{b_{i-1}}(i)<n_{b_{i-1}+1}(i)=n_{m_{i-1}+1}\qquad
  \text{and}\\
  n_{m_i}=n_{k(i)}(i)=b_i.
\end{gather*}
We deduce for $i\in\N$ that
\begin{align*}
  \Bigl\|\sum_{j=m_{i-1}+1}^{m_i}&\|P^E_{(n_{j-1},n_j]}(z_i)\|_Z\cdot
  v_j\Bigr \|_V\\
  =\ & \Bigl\|\sum_{j=b_{i-1}+1}^{k(i)} \|P^E_{(n_{j-1}(i),n_j(i)]}
    (z_i)\|_Z\cdot v_{j-(b_{i-1}-m_{i-1})}\Bigr \|_V\\
  =\ & \Bigl\|\sum_{j=a_i}^{k(i)} \|P^E_{(n_{j-1}(i),n_j(i)]}
    (z_i)\|_Z\cdot v_{j-(b_{i-1}-m_{i-1})}\Bigr \|_V\\
    & \text{(since $P^E_{(n_{j-1}(i),n_j(i)]}(z_i)=0$ for
    $j<a_i$)}\\[1ex]
  \geq\ & d\Bigl\|\sum_{j=a_i}^{k(i)} \|P^E_{(n_{j-1}(i),n_j(i)]}
    (z_i)\|_Z\cdot v_j\Bigr \|_V \geq cd\\
    & \text{(since $b_{i-1}-m_{i-1}\leq b_{i-1}$ and $L(b_{i-1})<a_i$).}
\end{align*}
Our claim now follows as in  the proof of Lemma \ref{L:2.1}.\end{proof}

\begin{lem}
  \label{L:2.5}
  Let $V$ and $U$ be  Banach spaces with normalized and
  $1$-unconditional bases $(v_i)$ and $(u_i)$, respectively, and
  assume that every subsequence of $(u_i)$ dominates every normalized
  block sequence of $(v_i)$ and that every subsequence of $(v_i)$ is
  dominated by every normalized block sequence of $(u_i)$. Let $Z$ be
  a Banach space with an FDD $(E_i)$.

  If $(E_i)$ satisfies  $U$-upper estimates in $Z$, then $(E_i)$ also
  satisfies $U$-upper estimates in $Z_V$.
\end{lem}
\begin{proof}
  It follows from the assumptions that for some constants $C_1,C_2$
  and $C_3$ in $[1,\infty)$ we have
    \begin{align}
      \label{E:2.5.1}
      &\Bigl\|\sum_{i=1}^\infty z_i\Bigr\|_Z \le C_1
      \Bigl\|\sum_{i=1}^\infty \|z_i\|_Z\cdot u_i \Bigr\|_U \text{ for all
	block sequences $(z_i)$ of $(E_i)$,}\\
      \label{E:2.5.2}
      &\text{subsequences of $(v_i)$ are $C_2$-dominated by normalized
	block sequences of $(u_i)$},\\
      \label{E:2.5.3}
      &\text{normalized block sequences of $(v_i)$ are $C_3$-dominated
	by subsequences of $(u_i)$.}
    \end{align}
    Let $K$ be the projection constant of $(E_i)$ in $Z$ and set
    $C=C_3+C_1C_2+2KC_2$. We  show that for any finite block sequence
    $(z_i)_{i=1}^m$ of $(E_i)$ and for any $k$ and $n_1<\ldots<n_k$
    in $\N$ we have (putting $z=\sum_{i=1}^m z_i$ and $n_0=0$)
    \begin{equation}
      \label{E:2.5.4}
      \Bigl\|\sum_{j=1}^k\|P^E_{(n_{j-1},n_j]}(z)\|_Z\cdot
      v_j\Bigl\|_V\le C\Bigl\|\sum_{i=1}^m\| z_i\|_{Z_V}\cdot
      u_i\Bigl\|_U.
    \end{equation}
    Taking then the supremum of the left side of~\eqref{E:2.5.4} over
    all choices of $k$ and $n_1<\ldots<n_k$ in $\N$, we obtain
    $$\Bigl\|\sum_{i=1}^m  z_i \Bigr\|_{Z_V}\le
    C\Bigl\|\sum_{i=1}^m\| z_i\|_{Z_V}\cdot u_i\Bigl\|_U,$$
    and thus that $(E_i)$ satisfies $C$-$U$-upper estimates in
    $Z_V$. Note that in proving~\eqref{E:2.5.4} we can of course
    assume that $n_k\leq\max\supp(z_m)$.
    
    For $i=1,2,\ldots,m$ put
    $$J_i=\big\{j\in\{1,2\ldots,k\}:\max\supp_E(z_{i-1})\le
    n_{j-1}<n_j\leq \max\supp_E(z_i)\big\}$$
    (with $\max\supp_E(z_0)=0$) and
    $J_0=\{1,2,\ldots,k\}\setminus\bigcup_{i=1}^m J_i$.
    
    For  $j=1,2,\ldots, k $ put
    $$I_j =\big\{i\in\{1,2,\ldots,m\}:
    n_{j-1}<\min\supp_E(z_i)<\max\supp_E(z_{i})\le n_j\big\}$$
    and $I_0=\{1,2,\ldots,m\}\setminus\bigcup_{j=1}^k I_j$.
    
    Firstly, we have
    \begin{align}\label{E:2.5.5}
      \Bigl\|\sum_{i=1}^m \sum_{j\in J_i}& \|P^E_{(n_{j-1},n_j]} (z_i)
	\|_{Z}\cdot v_j\Bigr\|_V\\
	&=\Bigl\|\sum_{i=1}^m b_i\Bigr\|_V\qquad\text{\big(where
	  $b_i=\sum_{j\in J_i} \|P^E_{(n_{j-1},n_j]} (z_i) \|_{Z}\cdot
	v_j$\quad for $1\leq i\le m$\big)} \notag\\
	  &\le C_3\Bigl\|\sum_{i=1}^m \|b_i\|_V\cdot u_i \Bigr\|_U
	  \notag\\
	  & \le C_3\Bigl\|\sum_{i=1}^m \|z_i\|_{Z_V}\cdot u_i
	  \Bigr\|_U\notag.
    \end{align}
    (Note that some (or all) of the $b_i$ may be zero, however the
    third line above is still valid using assumption~\eqref{E:2.5.3}.)
    Secondly,
    \begin{align}
      \label{E:2.5.6}
      \biggl\|\sum_{j\in J_0} \Bigl\|\sum_{i\in I_j}
      z_i\Bigr\|_Z&\cdot v_j\biggr\|_V\\
      \leq& C_1\Bigl\|\sum_{j\in J_0}\| b_j\|_U\cdot v_j\Bigr\|_V
      \notag\\
      & \text{(where $b_j=\sum_{i\in I_j}\| z_i\|_Z\cdot u_i$ for each
	$j\in J_0$, and we used~\eqref{E:2.5.1})}\notag\\
      \leq& C_1C_2\Bigl\|\sum_{j\in J_0}b_j\Bigr\|_U\notag\\
      \leq&C_1C_2\Bigl\|\sum_{i=1}^m\|z_i\|_{Z_V}\cdot
      u_i\Bigr\|_U,\notag
    \end{align}
    and, thirdly,
    \begin{align}\label{E:2.5.7}
      \Bigl\|\sum_{j\in J_0}
      \sum_{i\in I_0} \|&P^E_{(n_{j-1},n_j]} (z_i) \|_{Z}\cdot v_j\Bigr\|_V\\
	&\le \Bigl\|\sum_{j\in J_0}
	\|P^E_{(n_{j-1},n_j]} (z_{i^{(1)}_j}+z_{i^{(2)}_j})
	  \|_{Z}\cdot v_j\Bigr\|_V,\notag\\
	  \intertext{where the $z^{(1)}_{i_j}$'s and $z^{(2)}_{i_j}$'s
      are chosen as follows: since for every $j\in J_0$ the interval
      $(n_{j-1},n_j]$ intersects the support of at most two $z_i$'s
      with $i\in I_0$, we can choose, for $j\in J_0$,
      $i^{(1)}_j<i^{(2)}_j$ in $\{1,2,\ldots,m\}$ so that
      $i^{(2)}_j\le i^{(1)}_{j'}$ whenever $j<j'$ are in $J_0$ and so
      that the above inequality holds. Continuing~\eqref{E:2.5.7} we
      have}
	&\le \Bigl\|\sum_{j\in J_0}
	    \|P^E_{(n_{j-1},n_j]} (z_{i^{(1)}_j}) \|_{Z}\cdot v_j\Bigr\|_V+ \Bigl\|\sum_{j\in J_0}
	      \|P^E_{(n_{j-1},n_j]} (z_{i^{(2)}_j}) \|_{Z}\cdot v_j\Bigr\|_V\notag\\
		&\le KC_2 \Bigl\|\sum_{j\in J_0}
		\|z_{i^{(1)}_j} \|_{Z}\cdot u_{i^{(1)}_j}  \Bigr\|_U+KC_2 \Bigl\|\sum_{j\in J_0}
		\|z_{i^{(2)}_j} \|_{Z}\cdot u_{i^{(2)}_j}  \Bigr\|_U\notag\\
		&\le 2KC_2 \Bigl\|\sum_{i=1}^m \|z_i \|_{Z_V}\cdot u_i\Bigr\|_U\notag.
    \end{align}
    Finally, we deduce from~\eqref{E:2.5.5},~\eqref{E:2.5.6}
    and~\eqref{E:2.5.7} that
    \begin{align*}
      \Bigl\|\sum_{j=1}^k&\|P^E_{(n_{j-1},n_j]}(z)\|_Z\cdot v_j\Bigr\|_V  \\
	&\le \Bigl\|\sum_{i=1}^m\sum_{j\in J_i} \|P^E_{(n_{j-1},n_j]}(z_i)\|_Z\cdot v_j\Bigr\|_V  \\
	  &\qquad+\Bigl\|\sum_{j\in J_0}\|P^E_{(n_{j-1},n_j]}(z)\|_Z\cdot v_j\Bigr\|_V  \\
	    &\le \Bigl\|\sum_{i=1}^m\sum_{j\in J_i} \|P^E_{(n_{j-1},n_j]}(z_i)\|_Z\cdot v_j\Bigr\|_V  \\
	      &\qquad+\biggl\|\sum_{j\in J_0}\Bigl\|
	      \sum_{i\in I_j} z_i\Bigr\|_Z\cdot v_j\biggr\|_V  \\
		&\qquad+\Bigl\|\sum_{j\in J_0}\sum_{i\in I_0}\|P^E_{(n_{j-1},n_j]}(z_i)\|_Z\cdot v_j\Bigr\|_V  \\
		  &\le(C_3+C_1C_2+2KC_2)\Bigl\|\sum_{i=1}^m \|z_i \|_{Z_V}\cdot u_i\Bigr\|_U,
    \end{align*}
    which finishes the proof of   (\ref{E:2.5.4}).
\end{proof}

\section{Embedding Theorems}\label{S:3}

In this section we will prove and deduce some consequences of

\begin{thm}\label{T:3.1}
  Assume that $V$ is a Banach space with a normalized and
  1-unconditional basis $(v_i)$, and let $X$ be a  separable and
  reflexive space with $V$-lower tree estimates.
  \begin{enumerate}
  \item[a)]
    For every  reflexive space $Z$ with an FDD $E=(E_i)$ which
    contains $X$ there is a blocking $H=(H_i)$ of $(E_i)$ so that $X$
    naturally isomorphically embeds into $Z_V(H)$.
  \item[b)]
    There is a space $\Yt$ with a shrinking FDD $\Gt=(\Gt_i)$ so that
    $X$ is a quotient of $\Yt_V(\Gt)$.
  \end{enumerate}
\end{thm}
\begin{cor}\label{C:3.2}
  Assume that $V$ is a reflexive Banach space with a normalized and
  1-unconditional basis $(v_i)$ satisfying conditions  $(WLS)$ and
  $(SRS)$ as defined in the previous section and having the property
  that $(v_i)$ is dominated by every normalized block sequence of
  $(v_i)$. Let $X$ be a  separable and reflexive space with $V$-lower
  tree estimates.

  Then $X$ is a subspace of a reflexive space $Z$ with an FDD
  satisfying $V$-lower estimates and it is a quotient of a reflexive
  space $Y$ with an FDD satisfying $V$-lower estimates.
\end{cor}
\begin{remark}
  The assumption that $(v_i)$ is dominated by all its normalized block
  sequences implies that $(v_i)$ satisfies condition $(SRS)$.
\end{remark}
\begin{proof}
  By a theorem of Zippin~\cite{Z} we can embed $X$ into a reflexive
  space $W$ with an FDD $E=(E_i)$. Using Theorem~\ref{T:3.1}~(a) we
  can block $(E_i)$ into $F=(F_i)$ so that $X$ embeds into
  $Z=W_V(F)$. Theorem~\ref{T:3.1}~(b) provides a space $\Yt$ with a
  shrinking FDD $\Gt=(\Gt_i)$ so that $X$ is a quotient of
  $Y=\Yt_V(\Gt)$. By Corollary~\ref{C:2.4a} the spaces  $Z$ and  $Y$
  are reflexive.
 
  It follows from Lemma~\ref{L:2.3} that every normalized  block tree
  of $(F_i)$ in $Z$
  and of $(\Gt_i)$ in $Y$ has a branch which dominates some normalized
  block sequence of $(v_i)$ and thus $(v_i)$ itself. It follows that
  every normalized weakly null tree in $Z$ and in $Y$ has a branch
  which dominates $(v_i)$, and so, by  Proposition~\ref{P:1.2a}, $Z$
  and $Y$ satisfy $V$-lower  tree estimates. Finally, by
  Proposition~\ref{P:1.5} we can find blockings $G=(G_i)$ of $(F_i)$
  and $H=(H_i)$ of $(\Gt_i)$ so that $G$ satisfies $V$-lower
  estimates in $Z$, and $H$ satisfies $V$-lower  estimates in $Y$.
\end{proof}

From Corollary~\ref{C:3.2} and Proposition~\ref{P:1.3} we deduce in
certain instances the inverse implication of Proposition~\ref{P:1.4}.

\begin{cor}\label{C:3.2a}
  Assume that $V$ is a  reflexive  Banach space with  a
  1-unconditional and normalized basis $(v_i)$
  satisfying the conditions of Corollary~\ref{C:3.2}.
  
  If $X$ is a reflexive space which satisfies $V$-lower tree
  estimates, then $X^*$ satisfies $V^*$-upper tree estimates.
\end{cor}
\begin{thm}\label{T:3.1b}
  Let $V$ and $U$ be reflexive Banach spaces with normalized,
  $1$-unconditional bases $(v_i)$ and $(u_i)$, respectively, such that
  $(v_i)$ and $(u_i^*)$ both satisfy the conditions of
  Corollary~\ref{C:3.2}. Further assume that every subsequence of
  $(u_i)$ dominates every normalized block sequence of $(v_i)$ and
  that every normalized block sequence of $(u_i)$ dominates every
  subsequence of $(v_i)$.

  If $X$ is a separable, reflexive  Banach space which satisfies
  $(V,U)$-tree estimates, then $X$ can be embedded into a reflexive
  Banach space $Z$ with an FDD $(G_i)$, which satisfies
  $(V,U)$-estimates in $Z$.
\end{thm}
\begin{proof} By Proposition~\ref{P:1.4} $X^*$ satisfies $U^*$-lower
  tree estimates, and we can apply Corollary~\ref{C:3.2} to deduce
  that $X^*$ is the quotient of a  reflexive space $Y^*$ with an FDD
  $(E_i^*)$ ($Y^*$ being the dual of a space $Y$ with an FDD $(E_i)$)
  satisfying $U^*$-lower estimates in $Y^*$. Thus $X$ is a subspace of
  the reflexive space $Y$ having an FDD $(E_i)$ which, by
  Proposition~\ref{P:1.3}, satisfies $U$-upper estimates in $Y$.

  Theorem~\ref{T:3.1} part~(a) yields a blocking $F=(F_i)$ of $(E_i)$
  so that $X$ embeds into $Z=Y_V(F)$. As in the proof of
  Corollary~\ref{C:3.2} we can deduce from the assumptions that $Z$ is
  reflexive (Corollary~\ref{C:2.4a}), it satisfies $V$-lower tree
  estimates (Lemma~\ref{L:2.3} and Proposition~\ref{P:1.2a}) and that
  there is a blocking $G=(G_i)$ of $(F_i)$ such that $G$ satisfies
  $V$-lower estimates in $Z$ (Proposition~\ref{P:1.5}).

  To complete the proof note that the assumptions of Lemma~\ref{L:2.5}
  are satisfied, and so the FDD $(F_i)$, and hence also $(G_i)$,
  satisfies  $U$-upper estimates in $Z$.
\end{proof}

\begin{remark}
  Spaces $V$ which satisfy the assumptions of Corollary \ref{C:3.2}
  are the $\ell_p$ spaces, $1<p<\infty$, and convexified  Tsirelson
  spaces $T_{(p,\gamma)}$ (see \cite{CS}) for  $1\le p<\infty$ and
  $0<\gamma<1$. In section \ref{S:4} we will discuss more general
  versions of these spaces.
\end{remark}

The proof of Theorem \ref{T:3.1}  will follow along the lines
of the proof of Theorem 1.7  in \cite{OS2}, where the special
case $V=\ell_p$, for some $1<p<\infty$, was treated.

From Corollary 4.4 in \cite{OS1} we have
\begin{prop}
  \label{P:3.4}
  Let $X$ be a Banach space which is a subspace of a reflexive space $Z$
  with an FDD  $A=(A_i)$ having projection constant $K$. Let
  $\etab=(\eta_i)\subset (0,1)$ with $\eta_i\downarrow 0$. Then there exist
  positive integers $N_1 < N_2 < \ldots$ such that the following
  holds. For all $x\in S_X$ there exist $x_i \in X$ and $t_i\in
  (N_{i-1},N_i)$ ($i\in\N,\ N_0=0$) such that
  \begin{enumerate}
  \item[{a)}] $x = \sum\limits_{i=1}^\infty x_i$,
  \item[{b)}] for $i\in\N$ either $\|x_i\| <\eta_i$ or
    $\|P_{(t_{i-1},t_i)}^A x_i - x_i\| < \eta_i \|x_i\|$,
  \item[{c)}] $\|P_{(t_{i-1},t_i)}^A x-x_i\| < \eta_i$ for all $i\in\N$,
  \item[{d)}] $\|x_i\| < K+1$ for $i\in\N$, and
  \item[{e)}] $\|P_{t_i}^A x\| < \eta_i$ for $i\in\N$.
  \end{enumerate}
  Moreover, the above conditions hold if $(N_i)$ is replaced by any
  subsequence of $(N_i)$.
\end{prop}
Parts d) and e) were not explicitly stated in \cite{OS1} but follow from
the proof.

\begin{proof}[Proof of Theorem \ref{T:3.1} part (a)]
  Let $K$ be the projection constant of $E$ in $Z$, and assume that
  $X$ satisfies $C$-$V$-lower tree estimates.

  Using Proposition~\ref{P:3.3} as in the proof of
  Proposition~\ref{P:1.5}, we find a blocking $(F_i)$ of $(E_i)$ and
  a sequence $\deltab=(\delta_i)\subset (0,1)$, $\delta_i\downarrow 0$,
  such that any $\deltab$-skipped block sequence $(\xb_i)\subset S_X$
  of any blocking of $(F_i)$ in $Z$ $2C$-dominates $(v_i)$ \emph{and}
  $(v_{i+1})$. We can assume that
  $\Delta=\sum_i\delta_i<1$.

  It is easy to see that we can block $(F_i)$ into an FDD $G=(G_i)$ so
  that there exists $(e_n)\subset S_X$ with
  \begin{equation}
    \label{skipped-block}
    \| e_n-P^G_n(e_n)\|_Z<\delta_n\qquad\text{for all }
    n\in\N.
  \end{equation}

  Finally, we let $N_1<N_2<\ldots$ be a sequence of positive integers
  obtained by applying Proposition~\ref{P:3.4} with $(A_i)=(G_i)$ and
  $\etab=\deltab$.

  Now set $H_i=\bigoplus_{j=N_{i-1}+1}^{N_i}G_j$ for $i\in\N$, and let
  $H=(H_i)$. We consider the space $Z_V=Z_V(H)$, and claim that $X$
  naturally embeds into $Z_V$. In order to achieve that we prove that
  if $x\in S_X$, then
  \begin{equation}\label{E:3.1.4}
    \Big\|\sum_{i=1}^\infty \|P_i^H(x)\|_Z \cdot v_i\Big\|_V \leq
    24K^2C.
  \end{equation}
  Since the argument will also work for any further blocking of
  $(H_i)$ (by the ``moreover'' part of Proposition~\ref{P:3.4}) we
  obtain for all $x\in S_X$
  \begin{equation}\label{E:3.1.4a}
    \|x\|_{Z_V}=\sup_{0=k_0<k_1<\ldots<k_n} \Big\|\sum_{i=1}^n
    \|P^H_{(k_{i-1},k_i]}(x)\|_Z \cdot v_i\Big\|_V\leq 24K^2C.
  \end{equation}
  Let $x\in S_X$, and for each $i\in\N$ choose $x_i\in X$ and
  $t_i\in(N_{i-1},N_i)$ so that the properties (a)--(e) of
  Proposition~\ref{P:3.4} are satisfied with $(A_i)=(G_i)$ and 
  $\etab=\deltab$. 

  For each $i\in\N$ let $\xb_i=\frac{x_{i+1}}{\norm{x_{i+1}}}$ and
  $\alpha_i=\norm{x_{i+1}}$ if $\norm{x_{i+1}}\geq\delta_{i+1}$, and
  let $\xb_i=e_{N_i}$ and $\alpha_i=0$ if
  $\norm{x_{i+1}}<\delta_{i+1}$, where $(e_n)$ is a sequence that
  satisfies~\eqref{skipped-block}. Observe that $(\xb_i)$ is a
  $\deltab$-skipped block sequence of some blocking of $(G_i)$ (and
  hence of $(F_i)$) and as such it $2C$-dominates $(v_i)$ and
  $(v_{i+1})$. Using domination of $(v_i)$, we get
  \begin{align}
    \Bignorm{\sum_{i=1}^{\infty}x_i}_Z\geq &\Bignorm{\sum_{i=1}^\infty
      \alpha_i\xb_i}-\norm{x_1}-\Delta\\
    \geq &\frac{1}{2C}\Bignorm{\sum_{i=1}^\infty \alpha_i
      v_i}_V-(K+1)-\Delta\notag\\
    \geq &\frac{1}{2C}\Bignorm{\sum_{i=1}^\infty \norm{x_{i+1}}_Z\cdot
      v_i}_V- \frac{1}{2C}\Delta-(K+1)-\Delta\notag,
  \end{align}
  and thus
  \begin{equation}
    \label{dominate-vi}
    \Bignorm{\sum_{i=1}^\infty \norm{x_{i+1}}_Z\cdot v_i}_V\leq
    2C(K+2\Delta+2)\ .
  \end{equation}
  Since $(\xb_i)$ also $2C$-dominates $(v_{i+1})$, a similar
  calculation shows
  \begin{equation}
    \label{dominate-vi+1}
    \Bignorm{\sum_{i=1}^\infty \norm{x_{i+1}}_Z\cdot v_{i+1}}_V\leq
    2C(K+2\Delta+2)\ .
  \end{equation}
  Using properties~(c) and~(e) of Proposition~\ref{P:3.4} with $A=G$
  and $\etab=\deltab$, we have
  \[
  \norm{P^H_i(x)}_Z\leq K\norm{P^G_{(t_{i-1},t_{i+1})}(x)}_Z \leq
  K\big( \norm{x_i}_Z+\norm{x_{i+1}}_Z+3\delta_i\big)
  \]
  for each $i\in\N$. It follows, using~\eqref{dominate-vi}
  and~\eqref{dominate-vi+1}, that
  \begin{align*}
    \Bignorm{\sum_{i=1}^\infty \norm{P^H_i(x)}_Z\cdot v_i}_V &\leq
    \norm{P^H_1(x)}_Z+K \Bignorm{\sum_{i=2}^\infty
      \norm{x_i}_Z\cdot v_i}_V + K\Bignorm{\sum_{i=2}^\infty
      \norm{x_{i+1}}_Z\cdot v_i}_V + 3K\Delta\\
      &\leq 24K^2C\ .
  \end{align*}
\end{proof}
Before we prove part (b) of Theorem \ref{T:3.1} we need
a blocking result due to Johnson and Zippin.

\begin{prop}
  \label{P:3.6}
  \cite{JZ1}
  Let $T:Y\to Z$ be a bounded linear operator from a space $Y$ with a
  shrinking FDD  $(G_i)$ into a space $Z$ with an  $FDD$ $(H_i)$. Let
  $\vp_i\downarrow 0$. Then there exist blockings $E = (E_i)$ of $(G_i)$
  and $F = (F_i)$ of $(H_i)$ so that for all $m<n$ and $y\in
  S_{\bigoplus_{i\in (m,n)} E_i}$ we have $\|P_{[1,m)}^F Ty\| < \vp_m$
  and $\|P_{[n,\infty)}^F Ty\| <\vp_n$.
\end{prop}
%%%%%%%%%%%%%%%%%%%%%%%%%%%%%%%%%%%%%%%%%%%%%%%%%%%%%%%%%%%

\begin{proof}[Proof of Theorem \ref{T:3.1} part (b)]

  By Lemma 3.1 in~\cite{OS1} we can, after renorming $X$ if necessary,
  regard $X^*$ (isometrically) as a subspace of a reflexive space
  $Y^*$ (being the dual of a reflexive space $Y$) with bimonotone FDD
  $(E^*_i)$ such that $c_{00}(\oplus_{i=1}^\infty E^*_i)\cap X^*$ is
  dense in $X^*$. We have a natural quotient map $Q:Y\to X$. By a
  Theorem of Zippin~\cite{Z} we may regard $X$ (isometrically) as a
  subspace of a reflexive space $Z$ with an FDD $(F_i)$. Let $K$ be
  the projection constant of $(F_i)$ in $Z$, and choose $C>0$ such
  that $X$ satisfies $C$-$V$-lower tree estimates.

  Using Proposition~\ref{P:3.3} as in the proof of
  Proposition~\ref{P:1.5}, we find a sequence
  $\deltab=(\delta_i)\subset (0,1)$, $\delta_i\downarrow 0$, so that if
  $(x_i)\subset X$ is any $\deltab$-skipped block sequence of any
  blocking of $(F_i)$, then $(x_i)$ $2C$-dominates $(v_i)$, and
  moreover, using standard perturbation arguments and making $\deltab$
  smaller if necessary, we can assume that if $(z_i)\subset Z$
  satisfies $\|x_i-z_i\|<\delta_i$ for all $i\in\N$, then $(z_i)$ is a
  basic sequence equivalent to $(x_i)$ with projection constant at
  most $2K$. We also require that
  \begin{equation}
    \label{E:3.1.5}
    \Delta=\sum_{i=1}^{\infty}\delta_i<\frac17\ .
  \end{equation}
  Choose a sequence $\vpb=(\vp_i)\subset (0,1)$ with $\vp_i\downarrow 0$
  and
  \begin{equation}
    \label{E:3.3}
    3K(K+1)\sum_{j=i}^{\infty}\vp_j<\delta_i^2\qquad\text{for all
    }i\in\N\ .
  \end{equation}
  After blocking $(F_i)$ if necessary, we can assume
  that for any subsequent blocking $D$ of $F$ there is a sequence
  $(e_i)$ in $S_X$ such that
  \begin{equation}
    \label{skipped-block-again}
    \| e_i-P^D_i(e_i)\|_Z<\vp_i/2K\qquad\text{for all } i\in\N\ .
  \end{equation}
  By Proposition~\ref{P:3.6} we may assume, after further blocking our
  FDDs if necessary, that
  \begin{align}
    \label{E:3.1.6}
    &\text{for all $m<n$ and $y\in S_{\oplus_{i\in (m,n)} E_i}$ we
      have}\\
    &\qquad \|P^F_{[1,m)}\circ Q(y)\|<\vp_m\qquad\text{and}\qquad
      \|P^F_{[n,\infty)}\circ Q(y)\|<\vp_n\ ,\notag
  \end{align}
  and moreover the same holds if one passes to any blocking of
  $(E_i)$ and the corresponding blocking of $(F_i)$.

  For $i\in\N$ let $\Et_i$ be the quotient space of $E_i$ determined
  by $Q$, i.e.~if $y\in E_i$, then the norm of $\yt$, the equivalence
  class of $y$ in $E_i$, is given by $\trivert{\yt}=\|Q(y)\|$. Passing
  to a further blocking of $(E_i)$ (and the corresponding blocking of
  $(F_i)$), we may assume that $\Et_i\not=\{0\}$ for $i\in\N$. Given
  $y=\sum y_i \in \coo(\oplus_{i=1}^\infty E_i)$, $y_i\in E_i$ for
  $i\in\N$, we set $\yt=\sum \yt_i\in \coo(\oplus_{i=1}^\infty \Et_i)$
  and
  $$\trivert{\yt}=\max_{m\le n} \Big\|\sum_{i=m}^n
  Q(y_i)\Big\|=\max_{m\le n} \|Q\circ P^E_{[m,n]}(y)\|\ .$$
  We let $\Yt$ be the completion of $\coo(\oplus_{i=1}^\infty \Et_i)$
  with respect to $\trivert{\cdot}$. Since $(E_i)$ is a bimonotone FDD
  in $Y$, we have $\trivert{\yt}\leq\norm{y}$ for all $y\in
  \coo(\bigoplus_{i=1}^\infty E_i)$, and hence the map $y\mapsto \yt$
  extends to a norm one map from $Y$ to $\Yt$. By the definition of
  $\trivert{\cdot}$ we have $\norm{Qy}\leq \trivert{\yt}$ for any
  $y\in\coo(\oplus_{i=1}^\infty E_i)$. It follows that $\yt\to Q(y)$
  extends to a norm one map $\Qt\colon \Yt\to X$ with $\Qt(\yt)= Q(y)$
  for all $y\in Y$.
  
  In order to continue our proof we will need the following proposition
  from~\cite{OS2}.

  \begin{prop}
    \label{P:3.7}
    \cite[Proposition 2.6]{OS2}
    \begin{enumerate}
    \item[{a)}]
      $(\Et_i)$ is a bimonotone, shrinking FDD for $\Yt$.
    \item[{b)}]
      $\Qt$ is a quotient map from $\Yt $ onto $X$. More precisely if
      $x\in X$ and $y\in Y$ is such that $Q(y) = x$, $\|y\| =
      \|x\|$ and $y = \sum y_i$ with $y_i\in E_i$ for all $i\in\N$,
      then $\yt  = \sum \yt _i \in \Yt$, $\trivert{\yt } = \|y\|$ and
      $\Qt(\yt)  = x$.
    \item[{c)}] Let $(\yt _i)$ be a block sequence of $(\Et_i)$ in
      $B_{\Yt }$, and assume that $(\Qt(\yt _i))$ is a basic sequence
      with projection constant $\overline{K}$ and
      $a= \inf_i \|\Qt(\yt _i)\|>0$.
      Then for all  $(a_i)\in\coo$ we have
      \begin{equation*}
	\Bignorm{\sum a_i \Qt(\yt _i)}
	\le \Btrivert{\sum a_i\yt _i}
	\le \frac{3\overline{K}}{a} \Bignorm{\sum a_i \Qt(\yt _i)}\ .
      \end{equation*}
    \end{enumerate}
  \end{prop}

  To finish the proof of Theorem~\ref{T:3.1}~(b) it suffices to find a
  constant $L<\infty$ and a blocking $\Gt=(\Gt_i)$ of $(\Et_i)$ with
  the following property. For any $x\in S_X$ there exists a
  $\yt=\sum \yt_i\in \Yt$, $\yt_i\in \Gt_i$ for $i\in\N$, so that
  \begin{align}\label{E:3.7}
    &\|\Qt(\yt)-x\|<1/2\ ,\\
    &\label{E:3.8}\Big\|\sum_{j=1}^k\trivert{P^{\Gt}_{(n_{j-1},n_j]}(\yt)}
    \cdot v_j\Big\|_V\leq L\\
    &\ \text{for any choice of }k\ \text{and}\
    n_1<n_2<\ldots<n_k\text{ in }\N\ (n_0=0).\notag
  \end{align}
  Once this is accomplished, we consider the space
  $\Yt_V=\Yt_V(\Gt)$. Given $x=x_0\in S_X$, the property of $\Gt$
  allows us to recursively choose $x_n\in\frac1{2^n}B_X$ and $\yt_n\in
  \frac{L}{2^{n-1}}B_{\Yt_V}$, $n\in\N$, so that
  $x_n=x_{n-1}-\Qt(\yt_n)$ for all $n\in\N$. It follows that
  $\sum_{n=1}^\infty \yt_n $ converges in $\Yt_V$ with
  $\|\sum_{n=1}^\infty \yt_n\|\le 2L$ and $\Qt(\sum_{n=1}^\infty
  \yt_n)=x$. Thus $\Qt: \Yt_V\to X$  remains surjective, which
  finishes the proof.

  In order to show the existence of a suitable blocking $\Gt$ of $\Et$
  we need the following result from \cite{OS2}.
  
  \begin{lem}\label{L:3.8} \cite[Lemma  2.7]{OS2}
    Assume that~\eqref{E:3.1.6} holds for our original map $Q:Y\to
    X$. Then there exist integers $0= N_0 < N_1 <\dots$ so that if
    for each $i\in\N$ we define
    \begin{equation*}
      \begin{split}
	& C_i=\bigoplus_{j=N_{i-1}+1}^{N_i} E_j\ ,\qquad
	D_i=\bigoplus_{j=N_{i-1}+1}^{N_i} F_j,\\
	& L_i = \left\{ j\in \N : N_{i-1} < j\le  \frac{N_{i-1} +
	  N_i}2\right\}\ ,\\
	& R_i = \left\{ j\in \N : \frac{N_{i-1} +
	  N_i}2 <j\leq N_i \right\}\ ,\\
	& C_{i,L} = \bigoplus_{j\in L_i} E_j\qquad\text{and}\qquad
	C_{i,R} = \bigoplus_{j\in R_i} E_j \ ,
      \end{split}
    \end{equation*}
    then the following holds.
    Let $x\in S_X$, $0\leq m<n$ and $\vp>0$, and assume that
    $\norm{x-P^D_{(m,n)}(x)}<\vp$. Then there exists $y\in B_Y$ with
    $y\in C_{m,R}\oplus \Big(\bigoplus_{i\in (m,n)} C_i\Big)\oplus C_{n,L}$
    (where $C_{0,R}=\{0\}$) and $\norm{Qy - x}< K[2\vp +\vp_{m+1}]$
    (recall that $K$ is the projection constant of $(F_i)$ in $Z$).
  \end{lem}
  
  Let $(C_i)$ and $(D_i)$ be the blockings given by
  Lemma~\ref{L:3.8}. We now apply Proposition~\ref{P:3.4} with
  $(A_i)=(D_i)$ and $\etab=\vpb$ to obtain a sequence $N_1<N_2<\dots$
  in $\N$ so that the conclusions of the proposition are satisfied. We
  now come to our final blockings: for each $i\in\N$ set $G_i=\oplus
  _{j=N_{i-1}+1}^{N_i} C_j$ and let $H_i=\oplus _{j=N_{i-1}+1}^{N_i}
  D_j$ ($N_0=0$). Put $G=(G_i)$ and $H=(H_i)$. Let $\Gt=(\Gt_i)$ be
  the corresponding blocking of $(\Et_i)$.

  Fix a sequence $(e_i)$ in $S_X$ so that~\eqref{skipped-block-again}
  holds.
  Let $x\in S_X$. By the choice of $N_1, N_2, \dots$, for each
  $i\in\N$,  there are $x_i\in (K+1)B_X$ and $t_i\in(N_{i-1},N_i)$
  such that $x=\sum_{i=1}^\infty x_i$ and for all $i\in\N$ either
  $\|x_i\| <\vp_i$ or $\|P^D_{(t_{i-1},t_i)} x_i - x_i\| < \vp_i
  \|x_i\|$ ($t_0=0$). For each $i\in\N$ let
  $\xb_i=\frac{x_{i+1}}{\norm{x_{i+1}}}$ and $\alpha_i=\norm{x_{i+1}}$
  if $\norm{x_{i+1}}\geq \vp_{i+1}$, and let $\xb_i=e_{N_i}$ and
  $\alpha_i=0$ if $\norm{x_{i+1}}<\vp_{i+1}$.

  Since
  \begin{equation}
    \label{E:3.11}
    \norm{\xb_i-P^D_{(t_i,t_{i+1})}(\xb_i)}<\vp_{i+1}\qquad\text{for
      all}\ i\in\N\ ,
  \end{equation}
  there exists $(y_i)\subset B_Y$ with $y_i\in C_{t_i,R}\oplus
  \Big(\bigoplus_{j\in (t_i,t_{i+1})} C_j\Big)\oplus C_{t_{i+1},L}$
  and
  \begin{equation}
    \label{E:3.12}
    \norm{Q(y_i) -\xb_i}< 3K\vp_{i+1}\ ,\qquad i\in\N\ .
  \end{equation}
  Also, if
  $\norm{x_1}<\vp_1$, then set $y_0=0$, and if $\norm{x_1}\geq \vp_1$,
  then choose $y_0\in (K+1)B_Y$ such that $y_0\in \Big(\bigoplus_{j\in
  (0,t_1)} C_j\Big)\oplus C_{t_1,L}\subset G_1$ and $\norm{Q(y_0)
    -x_1}< 3K(K+1)\vp_1$.

  Set $\xb=x_1+\sum_{i=1}^{\infty} \alpha_i\xb_i$, and note that (this
  series converges and) by~\eqref{E:3.1.5} and~\eqref{E:3.3}
  \begin{equation}
    \label{E:3.4}
    \norm{x-\xb}\leq\sum_{i=2}^{\infty}\vp_i<\frac14\ .
  \end{equation}
  As a $\deltab$-skipped block sequence of a blocking of $(F_i)$
  (this follows from~\eqref{E:3.11} and~\eqref{E:3.3}),
  $(\xb_i)$ is a basic sequence that $2C$-dominates $(v_i)$. Since,
  by~\eqref{E:3.12},
  $\norm{\Qt(\yt_i) -\xb_i}< 3K\vp_{i+1}<\delta_i$ for all $i\in\N$,
  the sequence $\big(\Qt(\yt_i)\big)$ is also a basic sequence
  equivalent to
  $(\xb_i)$ with projection constant at most $2K$. Furthermore, we
  have $\inf_i\norm{\Qt(\yt_i)}\geq \inf_i
  \big(\norm{\xb_i}-\delta_i\big)>6/7$, and thus, by
  Proposition~\ref{P:3.7}~(c),
  \begin{equation}
    \label{E:3.13}
    \Bignorm{\sum a_i \Qt(\yt _i)} \leq \Btrivert{\sum a_i\yt _i}
	\leq 7K \Bignorm{\sum a_i \Qt(\yt _i)}\qquad\text{for all}\
	(a_i)\in\coo\ .
  \end{equation}
  Thus $(\yt_i)$ is a basic sequence equivalent to $(\xb_i)$ and, in
  particular, $\sum_{i=1}^{\infty} \alpha_i\yt_i$ converges. Putting
  $\yt=\yt_0+\sum_{i=1}^{\infty} \alpha_i\yt_i$ we have
  \begin{align*}
    \norm{\Qt\yt-\xb} &\leq \norm{\Qt\yt_0-x_1}+ \sum_{i=1}^{\infty}
    \abs{\alpha_i}\norm{\Qt\yt_i-\xb_i}\\
    &\leq 3K(K+1)\sum_{i=1}^{\infty}\vp_i<1/4\ ,
  \end{align*}
  and hence, by~\eqref{E:3.4}, $\norm{\Qt\yt-x}<1/2$, so we
  have~\eqref{E:3.7}.

  We now fix integers $k$ and $0=n_0<n_1<\dots<n_k$. For any $i\in\N$
  we have $P^{\Gt}_i(\yt)=P^{\Gt}_i(\yt_{i-1}+\yt_i)$. It follows that
  \begin{multline}
    \label{E:3.9}
    \Bignorm{\sum_{s=1}^k\trivert{P^{\Gt}_{(n_{s-1},n_s]}(\yt)} \cdot
    v_s}_V \leq \trivert{\yt_0}+
    \biggnorm{\sum_{s=1}^k\Btrivert{\sum _{%
	  \begin{subarray}{c}
	    i=n_{s-1}+1\\
	    i\geq 2
	\end{subarray}}^{n_s}
	\alpha_{i-1}\yt_{i-1}} \cdot v_s}_V\\
    + \biggnorm{\sum_{s=1}^k\Btrivert{\sum _{i=n_{s-1}+1}^{n_s}
	\alpha_i\yt_i} \cdot v_s}_V\ .
  \end{multline}
  We now show how to bound the third term of the right-hand side of
  the above inequality. A similar argument will give an estimate for
  the second term, and then~\eqref{E:3.8} will follow with
  $L=142K^2C$.

  For each $s=1,\dots,k$ let $\wt_s=\sum_{i=n_{s-1}+1}^{n_s}
  \alpha_i\yt_i$ and $b_s=\sum_{i=n_{s-1}+1}^{n_s}
  \alpha_i\xb_i$. Note that by~\eqref{E:3.11} and~\eqref{E:3.3}
  \begin{align}
    \label{E:3.10}
    \bignorm{P^D_{(t_{n_{s-1}+1},t_{n_s+1})}(b_s)-b_s} &\leq
    \sum_{i=n_{s-1}+1}^{n_s} \abs{\alpha_i}\cdot 2K\cdot
    \norm{P^D_{(t_i,t_{i+1})} \xb_i - \xb_i}\\
    &<2K(K+1)\sum_{i=n_{s-1}+1}^{n_s} \vp_{i+1}<\delta_s^2 &
    \text{($s=1,\dots,k$).}\notag
  \end{align}
  For $s>k$ set $n_s=n_k+(s-k)$. If $s>k$ or $\norm{b_s}<\delta_s$,
  then we let $\bb_s=\xb_{n_s}$ and $\beta_s=0$. If $s\leq k$ and
  $\norm{b_s}\geq\delta_s$, then we let $\bb_s=\frac{b_s}{\norm{b_s}}$
  and $\beta_s=\norm{b_s}$.

  It follows from~\eqref{E:3.10} and~\eqref{E:3.11} that $(\bb_s)$ is
  a $\deltab$-skipped block sequence of some blocking of $(F_i)$, and
  hence it is a basic sequence that $2C$-dominates
  $(v_i)$. From~\eqref{E:3.12} and~\eqref{E:3.13} we have
  \begin{align}
    \label{E:3.1}
    \norm{\Qt(\wt_s)-b_s} \leq& \sum_{i=n_{s-1}+1}^{n_s}
    \abs{\alpha_i}\norm{\Qt(\yt_i)-\xb_i}\\
    <&3K(K+1)\sum_{i=n_{s-1}+1}^{n_s}\vp_{i+1}
    <\delta_s\notag\\
    \intertext{and}
    \label{E:3.2}
    \trivert{\wt_s}\leq& 7K\norm{\Qt(\wt_s)} &\text{($s=1,\dots,k$).}
  \end{align}
  Putting $\gamma_s=\beta_s$ for $s=1,\dots,k$ and
  $\gamma_s=\alpha_{n_s}$ for $s>k$, we obtain the following.
  \begin{align*}
    \Bignorm{\sum _{s=1}^k \trivert{\wt_s}\cdot v_s}_V &\leq 7K
    \Bignorm{\sum _{s=1}^k \norm{\Qt(\wt_s)}\cdot v_s}_V
    &\text{(from~\eqref{E:3.2})}\\
    & \leq 7K \Bignorm{\sum _{s=1}^k \norm{b_s}\cdot
      v_s}_V+7K\Delta &\text{(from~\eqref{E:3.1})}\\
    &\leq 7K\Bignorm{\sum_{s=1}^{\infty} \gamma_sv_s}_V+14K\Delta
    &\text{(as $(v_s)$ is $1$-unconditional)}\\
    &\leq 7K\cdot2C\Bignorm{\sum_{s=1}^{\infty} \gamma_s\bb_s}
    +14K\Delta &\text{(since $(\bb_s)$ $2C$-dominates $(v_s)$)}\\
    &\leq 14KC\Bignorm{\sum_{i=1}^{\infty} \alpha_i\xb_i}+14KC\Delta
    +14K\Delta\\
    &\leq\text{\makebox[0pt][l]{$14KC\cdot (K+3)+14KC\Delta
      +14K\Delta<70K^2C\ ,$}}
  \end{align*}
  where the last line follows from~\eqref{E:3.4} and from
  $\norm{x_1}\leq K+1$. An almost identical calculation gives the same
  esimate for the second term in~\eqref{E:3.9}, and hence we
  obtain~\eqref{E:3.8} with $L=142K^2C$. This completes the proof of
  part~(b) of Theorem~\ref{T:3.1}.
\end{proof}

\section{Asymptotic estimates}\label{S:4}

Throughout this section we assume that $U$ and $V$ are Banach spaces
with normalized, $1$-subsymmetric bases $(u_i)$ and $(v_i)$,
respectively, i.e.~$(u_i)$ and $(v_i)$ are $1$-unconditional and
$1$-spreading ($\|\sum a_i v_i\|=\|\sum a_i v_{n_i}\|$ whenever
$(a_i)\in \coo$ and $n_1<n_2<\ldots$).
\begin{defin}
  Let $Z$ be a Banach space with an FDD $(E_i)$, and let $C\ge 1$. We
  say that $(E_i)$ {\em satisfies asymptotic $C$-$V$-lower estimates
    (in $Z$) }if for all $n\in\N$ we have:
  \begin{equation*}
    \Big\|\sum_{i=1}^n x_i\Big\|\ge C^{-1} \Big\|\sum_{i=1}^n
    \|x_i\|v_i\Big\|\quad\text{whenever $(x_i)_{i=1}^n$ is a block
      sequence of $(E_i)_{i=n}^\infty$}.
  \end{equation*}
  We say that $(E_i)$ {\em satisfies asymptotic $C$-$U$-upper
    estimates (in $Z$) }if for all $n\in\N$ we have:
  \begin{equation*}
    \Big\|\sum_{i=1}^n x_i\Big\|\le C \Big\|\sum_{i=1}^n
    \|x_i\|u_i\Big\|\quad\text{whenever $(x_i)_{i=1}^n$ is a block
      sequence of $(E_i)_{i=n}^\infty$}.
  \end{equation*}
  We say that $(E_i)$ {\em satisfies asymptotic $C$-$(V,U)$-estimates
    (in $Z$)} if it satisfies asymptotic $C$-$V$-lower and asymptotic
  $C$-$U$-upper estimates in $Z$.

  We say that $(E_i)$ {\em satisfies asymptotic $V$-lower estimates,
    $U$-upper estimates or $(V,U)$-estimates (in $Z$) }if there is a
  $C\geq 1$ so that $(E_i)$ satisfies asymptotic $C$-$V$-lower
  estimates, $C$-$U$-upper estimates or $C$-$(V,U)$-estimates in $Z$,
  respectively.
\end{defin}

As before we also introduce the coordinate-free version of asymptotic
lower and upper estimates, which can also be found (defined in  more
general situations) in~\cite{MT} and~\cite{MMT}.

\begin{defin}
  We say that a reflexive space $X$ {\em satisfies asymptotic
    $C$-$V$-lower tree estimates }or {\em asymptotic $C$-$U$-upper
    tree estimates }if, for every $k$, every normalized weakly null
  tree of length $k$ in $X$ has a branch which $C$-dominates
  $(v_i)_{i=1}^k$ or is $C$-dominated by $(u_i)_{i=1}^k$,
  respectively. We say $X$ {\em satisfies asymptotic $C$-$(V,U)$-tree
    estimates }if it satisfies asymptotic $C$-$V$-lower tree estimates
  and asymptotic $C$-$U$-upper tree estimates.

  We will say  that $X$ {\em satisfies  asymptotic $V$-lower tree,
    $U$-upper tree or $(V,U)$-tree estimates }if there is a $C\ge 1$
  so that $X$ satisfies asymptotic $C$-$V$-lower tree, $C$-$U$-upper
  tree or $C$-$(V,U)$-tree estimates, respectively.
\end{defin}
The following dualities can be shown  as in Propositions~\ref{P:1.3}
and~\ref{P:1.4}.
\begin{prop}\label{P:4.0}
  Assume that $V$ is a space with a normalized, 1-subsymmetric basis
  $(v_i)$.
  \begin{enumerate}
  \item[a)]
    For a space $Z$ with an FDD $(E_i)$ the following statements are
    equivalent.
    \begin{enumerate}
    \item[1)]
      $(E_i)$ satisfies asymptotic $V$-lower estimates in $Z$.
    \item[2)]
      $(E^*_i)$ satisfies asymptotic $\Vs$-upper estimates in  $\Zs$.
    \end{enumerate}
  \item[b)]
    If $X$ is a reflexive space which satisfies asymptotic
    $C$-$V$-upper tree estimates for some $C\ge 1$, then, for any
    $\vp>0$, $X^*$ satisfies asymptotic $(2C+\vp)$-$\Vs$-lower tree
    estimates.
\end{enumerate}
\end{prop}

For a space $V$ with a normalized, $1$-subsymmetric basis $(v_i)$ and
for $0<\gamma<1$ we will introduce the {\em Tsirelson space
  $T(V,\gamma)$ associated to $V$ and $\gamma$} as follows. It is the
space defined as the completion of $\coo$ under the norm
$\norm{\cdot}_{T(V,\gamma)}$, where
$$\norm{x}_{T(V,\gamma)}=\max_{\ell\in\N_0}\norm{x}_{\ell,T(V,\gamma)}
\qquad\text{for all }x\in \coo\ ,$$
and the norms $\norm{\cdot}_{\ell,T(V,\gamma)},\ \ell\in\N$, on $\coo$
are defined recursively as follows. For $x=(x_i)\in\coo$ we put
\begin{align*}
  &\norm{x}_{0, T(V,\gamma)}=\norm{x}_\infty=\max_{i\in\N}|x_i|\\
  \intertext{and, assuming $\norm{\cdot}_{\ell, T(V,\gamma)}$ has been
    defined, we put}
  &\norm{x}_{\ell+1, T(V,\gamma)}=\norm{x}_{\ell, T(V,\gamma)}\vee
  \max_{n\in\N,\ n\leq A_1<A_2<\dots<A_n} \gamma\Bignorm{\sum_{i=1}^n
    \norm{P_{A_i}(x)}_{\ell, T(V,\gamma)}\cdot v_i}_V\ ,
\end{align*}
where for $A,B\subset \N$ and $n\in\N$, $n\leq A$ means that $n\leq a$
for all $a\in A$, and $A<B$ means that $a<b$ for all $a\in A$ and
$b\in B$. $P_A$, for $A\subset \N$, denotes the projection
$\sum a_i e_i\mapsto \sum_{i\in A} a_i e_i$ on $\coo$.

As in the case of $V=\ell_1$, which yields the standard Tsirelson
space (cf.~\cite{CS}), it is easy to see that
$\norm{\cdot}_{T(V,\gamma)}$ satisfies the following implicit
equation.
\begin{equation}\label{E:4.1}
  \norm{x}_{T(V,\gamma)}=\norm{x}_\infty\vee \sup_{n\ge2,\ n\le
    A_1<A_2<\dots<A_n} \gamma \Bignorm{\sum_{i=1}^n
    \norm{P_{A_i}(x)}_{T(V,\gamma)}\cdot v_i}_V
\end{equation}
whenever $x\in T(V,\gamma)$.
\begin{prop}
  \label{P:4.3}
  Let $V$ be a Banach space with a normalized, 1-subsymmetric basis
  $(v_i)$. Let $X$ be a Banach space with a normalized basis $(e_i)$
  satisfying asymptotic $C$-$V$-lower estimates for some $C\geq 1$,
  i.e.
  \begin{equation}\label{E:4.3.1}
    \Bignorm{\sum_{i=1}^n a_i x_i}\geq C^{-1}\Bignorm{\sum_{i=1}^n
      a_iv_i}_V
  \end{equation}
  whenever $n\in\N,\ (x_i)_{i=1}^n$ is a normalized block sequence of
  $(e_i)_{i=n}^\infty$ in $X$ and $(a_i)_{i=1}^n\subset \R.$ Let $K$
  be the projection constant of $(e_i)$ in $X$. Then $(e_i)$
  $K$-dominates the unit vector basis $(t_i)$ of $T(V,\gamma)$
  whenever $0<\gamma<1$ and $\gamma\leq (KC)^{-1}$.
\end{prop}
\begin{proof} There is an equivalent norm $\trivert{\cdot}$ on $X$
  with respect to which $(e_i)$ is bimonotone and which satisfies
  $\norm{x}\leq\trivert{x}\leq K\norm{x}$ for all $x\in X$. In
  $\trivert{\cdot}$ the basis $(e_i)$ satisfies asymptotic
  $(KC)$-$V$-lower estimates. We can easily show by induction on
  $l\in\N_0$ that
  $$\Btrivert{\sum a_i e_i}\geq\Bignorm{\sum a_i
    t_i}_{\ell,T(V,\gamma)}\qquad\text{whenever } (a_i)\in\coo\ ,$$
  which proves the proposition.
\end{proof}

\begin{prop}\label{P:4.1}
  Let $V$ be a Banach space with a normalized, 1-subsymmetric basis
  $(v_i)$, and let $0<\gamma<1$.
  \begin{enumerate}
  \item[a)]
    The  unit vector basis $(t_i)$ in $T(V,\gamma)$ is a normalized,
    1-unconditional basis of $T(V,\gamma)$.
  \item[b)]
    If $(v_i)$ is not equivalent to the unit vector basis of $\co$
    (i.e.~$\bignorm{\sum_{i=1}^n v_i}\to\infty$ as
    $n\to\infty$), then $T(V,\gamma)$ is reflexive.
  \item[c)]
    $(t_i)$ is 1-dominated by every normalized block sequence of
    $(t_i)$ and, in particular (as subsequences are normalized block
    sequences), $T(V,\gamma)$ satisfies (SRS) with $c=1$.
  \item[d)]
    If $\gamma<1/4$, then $T(V,\gamma)$ satisfies (WLS) for any
    $d<1$.
  \end{enumerate}
\end{prop}
\begin{proof}
  (a) is trivial and (b) follows from the fact that $T(V,\gamma)$ does
  not contain $\ell_1$ or $\co$, which can be shown as in the case
  $V=\ell_1$ (cf.~\cite{CS}), and from the theorem of James
  (cf.~\cite{LT}) that states that a space with an unconditional
  basis is reflexive if and only if it does not contain copies of
  $\co$ or $\ell_1$.

  For a normalized block sequence $(x_i)$ of $(t_i)$ it follows
  from~\eqref{E:4.1} that
  $$ \Bignorm{\sum a_i x_i}\geq \max_{n,\ n\leq k_0<k_1<\dots<k_n}
  \gamma\biggnorm{\sum_{j=1}^n \Bignorm{\sum_{i=k_{j-1}}^{k_j-1} a_i
      x_i}\cdot v_j}_V$$
  whenever $(a_i)\in\coo$. Thus $(x_i)$ is a normalized, bimonotone
  basic sequence satisfying asymptotic $C$-$V$-lower estimates for
  $C=\gamma^{-1}$, and hence (c) follows from
  Proposition~\ref{P:4.3}.

  Claim (d) follows from the following lemma.
\end{proof}
\begin{lem}\label{L:4.2}
  For $m\in\N$, let $S^m\colon\coo\to\coo$ be the backward shift
  by $m$ coordinates, i.e.
  $$S^m\Big(\sum_{i=1}^\infty a_ie_i\Big)=\sum_{i=1}^\infty
  a_{i+m}\,e_i\qquad\text{for all } (a_i)\in\coo\ .$$
  Assume that $0<\gamma<1/4$ and that $m<n$ in $\N$ satisfy
  $$\frac{m}{n-m}\cdot \frac1{1-4\gamma}<1\ .$$
  Then for any $x\in \coo$ with $n\le \min\supp(x)$ we have
  $$\norm{S^m(x)}_{T(V,\gamma)}\geq\Big( 1-\frac{m}{n-m}
  \cdot\frac1{1-4\gamma}\Big)\norm{x}_{T(V,\gamma)}\ .$$
\end{lem}
\begin{proof}
  Set  $\norm{\cdot}=\norm{\cdot}_{T(V,\gamma)}$ and
  $\norm{\cdot}_\ell=\norm{\cdot}_{\ell,T(V,\gamma)}$ for
  $\ell\in\N$. Fix $m\in\N$. Given $\ell,n\in\N$ with $m<n$ and
  $\frac{m}{n-m}\cdot\frac1{1-4\gamma}<1$, we put
  $$\Delta(0,n)=0\qquad\text{and}\qquad\Delta(\ell,n)=
  \frac{m}{n-m}\sum_{i=0}^{\ell-1} (4\gamma)^i\ .$$
  We will show by induction on $\ell\in\N_0$ that
  \begin{multline}\label{E:4.2.1}\allowdisplaybreaks[1]
    \|S^m(x)\|_\ell\geq\big(1-\Delta(\ell,n)\big)\cdot\norm{x}_\ell\\
    \text{whenever } x\in\coo,\ m<n\leq\min\supp(x)\ \text{and}\
    \frac{m}{n-m}\cdot \frac1{1-4\gamma}<1\ .
  \end{multline}
  Clearly,~\eqref{E:4.2.1} yields the lemma, and  the inequality is
  trivially true for $\ell=0$. Assume that~\eqref{E:4.2.1} is true for
  some $\ell\in\N_0$. Fix $n\in\N$ and $x\in\coo$ such that
  $m<n\leq\min\supp(x)$, $\frac{m}{n-m}\cdot \frac1{1-4\gamma}<1$
  and $\|x\|_{\ell+1}=1$.

  If $\norm{x}_{\ell+1}=\norm{x}_\ell$, then~\eqref{E:4.2.1} follows for
  $\ell+1$ from the induction hypothesis. Otherwise we can find $k\in\N$,
  $k\ge n$, and sets $ k\leq A_1<A_2<\dots<A_k$ so that
  $$\|x\|_{\ell+1}=\gamma\Big\|\sum_{i=1}^k \|P_{A_i}(x)\|_\ell \cdot
  v_i\Big\|_V=1$$
  with  $\|P_{A_i}(x)\|_\ell\le1$ for $i=1,2,\ldots,k$ (otherwise
  $\|x\|_\ell\ge \|P_{A_i}(x)\|_\ell > 1= \|x\|_{\ell+1} $).

  Let $I\subset \{1,2,\ldots,k\}$, $\#I=k-m$, so that
  $$\Big\|\sum_{i\in I} \|P_{A_i}(x)\|_\ell\cdot v_i\Big\|_V=
  \max_{\tilde I\subset\{1,2\ldots,k\},\ \#\tilde I=k-m}
  \Big\|\sum_{i\in\tilde I} \|P_{A_i}(x)\|_\ell\cdot v_i\Big\|_V\ .$$
  Then it follows that
  $$\gamma\Big\|\sum_{i\in I} \|P_{A_i}(x)\|_\ell\cdot v_i\Big\|_V\geq
  \frac{k-m}{k}\geq \frac{n-m}{n}\ .$$
  Indeed
  \begin{equation*}
    \left(\begin{matrix} k\\ k-m\end{matrix}\right)\cdot
      \Big\|\sum_{i\in I} \|P_{A_i}(x)\|_\ell\cdot v_i\Big\|_V
      \ge\Big\|\sum_{%
	\begin{subarray}{c}
	  \tilde I\subset\{1,2\ldots ,k\},\\
	  \#\tilde I=k-m
      \end{subarray}}%
      \sum_{i\in\tilde  I} \|P_{A_i}(x)\|_\ell\cdot v_i\Big\|_V=
      \frac1\gamma\left(\begin{matrix} k-1\\
	k-m-1\end{matrix}\right)\ ,
  \end{equation*}
  and thus
  \begin{equation*}
    \gamma\Big\|\sum_{i\in I} \|P_{A_i}(x)\|_\ell\cdot v_i\Big\|_V\ge
    \left(\begin{matrix} k-1\\ k-m-1\end{matrix}\right)
      \Big/\left(\begin{matrix} k\\ k-m\end{matrix}\right)
	=\frac{k-m}{k}\ge \frac{n-m}{n}\ .
  \end{equation*}
  Put
  \begin{align*}
    &I_1=\{i\in I : \#A_i\le\min A_i-m\}\qquad\text{and}\\
    &I_2=I\setminus I_1=\{i\in I : \#A_i >\min
    A_i-m\}\ .
  \end{align*}
  For $i\in I_1$ we have
  $$\|S^{m}(P_{A_i}(x))\|_\ell=\|P_{A_i}(x)\|_\ell\ .$$
  Let $k'\leq k-m$ and $i_1<i_2<\dots<i_{k'}$ be such that
  $I_2=\{i_1,i_2,\ldots,i_{k'}\}$. We note that
  \begin{align*}
    \min(A_{i_1})&\ge k\ge n\ ,\\
    \min(A_{i_2})&\ge \min(A_{i_1})+\#A_{i_1}>
    2\cdot\min(A_{i_1})-m\ge 2n-m\ ,\\
    \intertext{and more generally, by induction, for each
      $s=2,\dots,k'$ we have}
    \min(A_{i_s})&\geq\min(A_{i_{s-1}})+\#A_{i_{s-1}}>
    2\cdot\min(A_{i_{s-1}})-m \\
    & \geq 2\cdot\big(2^{s-2} n -(2^{s-2}-1)m\big)-m
    =2^{s-1}n-(2^{s-1}-1)m\ .
  \end{align*}
  We deduce that
  \begin{align*}
    \|S^m&(x)\|_{\ell+1}\\
    &\geq \gamma\Big\|\sum_{i\in I} \|S^m(P_{A_i}(x))\|_\ell \cdot
    v_i\Big\|_V\\
    &\text{(note that $\#I=k-m\le \min A_1 -m$)}\\
    &\ge \gamma\Big\| \sum_{i\in I_1} \|P_{A_i}(x)\|_\ell\cdot v_i
    \!+\!\sum_{s=1}^{k'}\Big(1-
    \Delta\big(\ell,2^{s-1}n-(2^{s-1}-1)m\big)\Big)
    \|P_{A_{i_s}}(x)\|_\ell\cdot v_{i_s}\Big\|_V\\
    &\text{(by the induction hypothesis and since $(v_i)$ is
      1-unconditional)}\\
    &\ge \gamma\Big\| \sum_{i\in I} \|P_{A_i}(x)\|_\ell\cdot v_i\Big\|_V
    -\gamma \sum_{s=1}^{k'}
    \Delta\big(\ell,2^{s-1}n-(2^{s-1}-1)m\big)\\
    &\text{(since $\|P_{A_i}(x)\|_\ell\le 1$ for $i=1,2,\dots,k$)}\\
    &\ge \frac{n-m}{n} -\gamma \sum_{s=1}^{k'} \sum_{t=0}^{\ell-1}
    (4\gamma)^t \frac{m}{2^{s-1}(n-m)}\\
    &= 1- \frac{m}{n}-\gamma\sum_{t=0}^{\ell-1}  (4\gamma)^t
    \frac{m}{(n-m)}\sum_{s=1}^{k'} 2^{1-s}\\
    &\ge 1-\frac{m}{n-m}-4\gamma\sum_{t=0}^{\ell-1}  (4\gamma)^t
    \frac{m}{(n-m)} =1-\Delta(\ell+1,n)\ ,
  \end{align*}
  which finishes the proof of our induction step.
\end{proof}

\begin{prop}\label{P:4.4}
 Let $V$ be a space with a normalized and 1-subsymmetric basis
 $(v_i)$. For a separable and reflexive Banach space $X$ the following
 are equivalent.
 \begin{enumerate}
 \item[a)]
   $X$ satisfies asymptotic $V$-lower tree estimates.
 \item[b)]
   There exists $\gamma\in(0,1)$ such that $X$ satisfies
   $T(V,\gamma)$-lower tree estimates.
 \end{enumerate}
\end{prop}
\begin{proof}
  Using a result of Zippin~\cite{Z}, and after renorming $X$ if
  necessary, we can assume that $X$ is (isometrically) a subspace of a
  reflexive space  $Z$ with a bimonotone FDD $(E_i)$.
  
  \noindent
  ``(a)$\Rightarrow$(b)'' Let $C\ge 1$ so that $X$ satisfies
  asymptotic $C$-$V$-lower tree estimates. For $k\in\N$ set
  \[
  \cA^{(k)}=\Big\{(x_i)\in S_X^\omega \,:\,\Big\|\sum_{i=1}^k a_i
  v_i\Big\|_V\le C \Big\|\sum_{i=1}^k a_i x_i\Big\|_X \text{ for all
  }(a_i)_{i=1}^k\subset \R\Big\}\ .
  \]
  Let $\vp>0$ be small enough so that for all $k\in\N$
  \[
  \overline{\cA^{(k)}_\vp}\subset
  \Big\{(x_i)\in S_X^\omega \,:\,\Big\|\sum_{i=1}^k a_i v_i\Big\|_V\le
  2C \Big\|\sum_{i=1}^k a_i x_i\Big\|_X \text{ for all
  }(a_i)_{i=1}^k\subset \R\Big\}\ .
  \]
  We let $(E^{(0)}_i)=(E_i)$ and apply Proposition~\ref{P:3.3} to
  $\cA^{(k)}$ successively for each $k\in\N$ to obtain decreasing
  null sequences $\deltab^{(k)}=(\delta^{(k)}_i)\subset(0,1)$ and
  blockings $(E^{(k)}_i)$ of $(E^{(k-1)}_i)$ so that if $(x_i)\subset
  S_X$ is a $\deltab^{(k)}$-skipped block sequence of $(E^{(k)}_i)$,
  then $(x_i)$ lies in $\overline{\cA^{(k)}_\vp}$.

  Let $(G_i)$ be a blocking of $(E_i)$ such that $(G_i)_{i=k}^\infty$
  is a blocking of $(E^{(k)}_i)_{i=k}^\infty$ for all $k\in\N$. Then
  choose $\deltab=(\delta_i)\subset(0,1)$, $\delta_i\downarrow 0$,
  such that if $(x_i)\subset S_X$ is a $\deltab$-skipped block
  sequence of $(G_i)$, then $(x_i)$ is a basic sequence with
  projection constant at most $2$ and for any $k\in\N$ any normalized
  block sequence $(z_i)_{i=1}^k$ of $(x_i)_{i=k}^\infty$ is a
  $\deltab^{(k)}$-skipped block sequence of $(E^{(k)}_i)$. It follows
  that any $\deltab$-skipped block sequence of $(G_i)$ is a basic
  sequence with projection constant at most~$2$ satisfying
  asymptotic $(2C)$-$V$-lower estimates, and hence, by
  Proposition~\ref{P:4.3}, it $2$-dominates the unit vector basis
  $(t_i)$ of $T(V,\gamma)$ for $\gamma=(4C)^{-1}$. Thus (b) follows.

  \noindent
  ``(b)$\Rightarrow$(a)'' This follows from Proposition~\ref{tree-est}
  and the fact that, by~\eqref{E:4.1}, $(t_i)_{i=k}^{2k-1}$
  $C$-dominates $(v_i)_{i=1}^k$ with $C=\gamma^{-1}$, where $(t_i)$ is
  the unit vector basis of $T(V,\gamma)$.
\end{proof}
We are now ready to state the main results of this section.
\begin{thm}\label{T:4.5}
  Assume that $V$ is a space  with a normalized and 1-subsymmetric
  basis  $(v_i)$ and that $X$ is a separable, reflexive space with
  asymptotic $V$-lower tree estimates. Then $X$ can be embedded in a
  reflexive space $Z$ with an FDD satisfying asymptotic $V$-lower
  estimates, and $X$ is isomorphic to a quotient of a reflexive space
  $Y$ with an FDD satisfying  asymptotic  $V$-lower estimates.
\end{thm}
\begin{proof}
  If $(v_i)$ is equivalent to the unit vector basis of $\co$, then we
  simply use the theorem of Zippin~\cite{Z} to embed $X$ and $X^*$
  into reflexive spaces $Z$ and $Y^*$, respectively, with FDDs. The
  result then follows, since any FDD satisfies (asymptotic)
  $\co$-lower estimates.

  Assume now that $(v_i)$ is not equivalent to the unit vector basis
  of $\co$. By Proposition~\ref{P:4.4} $X$ satisfies
  $T(V,\gamma)$-lower tree estimates for some $\gamma\in (0,1)$. We
  may clearly assume that $\gamma<1/4$, and then, by
  Proposition~\ref{P:4.1}, $T(V,\gamma)$ is a reflexive space whose
  unit vector basis $(t_i)$ is a normalized and $1$-unconditional
  basis satisfying $(SRS)$ and $(WLS)$ and dominated by its normalized
  block sequences. Hence, by Corollary~\ref{C:3.2}, $X$ embeds into a
  reflexive space $Z$ with an FDD satisfying $T(V,\gamma)$-lower
  estimates, and $X$ is isomorphic to a quotient of a reflexive space
  $Y$ with an FDD satisfying $T(V,\gamma)$-lower estimates. It is
  clear (e.g.~from~\eqref{E:4.1}) that the FDDs of $Z$ and $Y$ satisfy
  asymptotic $V$-lower esimates.
\end{proof}

\begin{thm}\label{T:4.6}
  Let $V$ and $U$ be Banach spaces with normalized, $1$-subsymmetric
 bases $(v_i)$ and $(u_i)$, respectively. Assume that, for any
 $\gamma\in(0,1/4)$, every normalized block sequence of the unit
 vector basis of $\big(T(\Us,\gamma)\big)^*$ dominates every
 normalized block sequence of the unit vector basis of
 $T(V,\gamma)$. For a separable and reflexive Banach space $X$ the
 following are equivalent.
  \begin{enumerate}
  \item[a)]
    $X$ satisfies asymptotic $(V,U)$-tree estimates.
  \item[b)]
    $X$ can be embedded in a reflexive space $Z$ with an FDD $(E_i)$
    satisfying asymptotic $(V,U)$-estimates.
  \item[c)]
    $X$  is the quotient of a reflexive space $Z$ with an FDD $(E_i)$
    satisfying asymptotic  $(V,U)$-estimates.
  \item[d)]
    $X^*$ satisfies asymptotic $(\Us,\Vs)$-tree estimates.
  \end{enumerate}
\end{thm}
\begin{remark}
  The conditions of the theorem are satisfied by certain pairs of
  $\ell_p$ spaces. This will be spelt out in Corollary~\ref{C:4.7}
  below. Note that $T(\co,\gamma)$ is just the space $\co$ for any
  $\gamma\in(0,1)$. So if one of the bases $(v_i)$ and $(u_i^*)$ of
  $V$ and $\Us$, respectively, are equivalent to the unit vector basis
  of $\co$, then the assumptions on the spaces $T(V,\gamma)$ and
  $\big(T(\Us,\gamma)\big)^*$ are automatically satisfied. In this
  case the theorem is really a statement about one-sided estimates
  since every reflexive space satisfies $(\co,\ell_1)$-tree estimates,
  and any FDD satisfies $(\co,\ell_1)$-estimates.
\end{remark}
\begin{proof}
  First note that the implication ``(b)$\Rightarrow$(a)'' is clear
  and, by Proposition~\ref{P:4.0}~(a), that the implication
  ``(c)$\Rightarrow$(d)'' is just another instance of
  ``(b)$\Rightarrow$(a)''. Also, the implications
  ``(a)$\Rightarrow$(d)''  and ``(d)$\Rightarrow$(a)'' are
  equivalent since the pair $(\Us,\Vs)$ satisfies the same assumptions
  as the pair $(V,U)$. Thus we will have completed the proof once we
  show how to deduce~(b) and~(c) from~(a).

  The assumption that~(a) holds splits into the following two
  conditions by Proposition~\ref{P:4.0}~(b): $X$ satisfies asymptotic
  $V$-lower tree estimates, and $X^*$ satisfies asymptotic $\Us$-lower
  tree estimates. If $(u_i^*)$ is equivalent to the unit vector basis
  of $\co$, then the second condition is redundant (cf.~the remark
  preceding this proof), and~(b) and~(c) follow from an application of
  Theorem~\ref{T:4.5} to the pair $(V,X)$. Similarly, if $(v_i)$ is
  equivalent to the unit vector basis of $\co$, then the first
  condition is redundant, and~(b) and~(c) follow from an application
  of Theorem~\ref{T:4.5} to the pair $(\Us,X^*)$ followed by an
  application of Theorem~\ref{P:4.0}~(a). In general, by
  Proposition~\ref{P:4.4}, it follows from the two conditions that,
  for some $\gamma\in(0,1/4)$, $X$ satisfies 
  $T(V,\gamma)$-lower tree estimates, and $X^*$ satisfies 
  $T(\Us,\gamma)$-lower tree estimates.

  We now continue the proof under the assumption that neither $(v_i)$
  nor $(u_i^*)$ is equivalent to the unit vector basis of $\co$. Then
  it follows from Proposition~\ref{P:4.1} that $T(V,\gamma)$ and
  $T(\Us,\gamma)$ are reflexive spaces each having a normalized,
  $1$-unconditional basis satisfying $(SRS)$ and $(WLS)$ and dominated
  by all its normalized block sequences. Thus we are in the situation
  of Corollary~\ref{C:3.2a}, two applications of which give that 
  $X$ satisfies $\big(T(V,\gamma),T(\Us,\gamma)^*\big)$-tree
  estimates, and $X^*$ satisfies
  $\big(T(\Us,\gamma),T(V,\gamma)^*\big)$-tree estimates.

  We now complete the proof by applying Theorem~\ref{T:3.1b} first to
  $\big(T(V,\gamma),T(\Us,\gamma)^*\big)$ and $X$ to deduce~(b) and
  then to $\big(T(\Us,\gamma),T(V,\gamma)^*\big)$ and $X^*$ to
  obtain~(c) (the second application of Theorem~\ref{T:3.1b} is
  followed by an application of Proposition~\ref{P:4.0}~(a)). Note
  that we are assuming that every normalized block sequence of the
  unit vector basis of $\big(T(\Us,\gamma)\big)^*$ dominates every
  normalized block sequence of the unit vector basis of $T(V,\gamma)$
  (which then implies the same for the spaces $T(V,\gamma)^*$ and
  $(T(\Us,\gamma)$), so the conditions of Theorem~\ref{T:3.1b} are
  indeed satisfied.
\end{proof}

Let $0<\gamma<1$ and $1\leq p<\infty$. We shall write $T_{p,\gamma}$
for the Tsirelson space $T(\ell_p,\gamma)$ associated to $\ell_p$ and
$\gamma$.
\begin{cor}\label{C:4.7}
  Let $1\le q\le p \le\infty$, and let $\frac1p+\frac1{p'}=1$ and
  $\frac1q+\frac1{q'}=1$. For a separable and reflexive Banach space
  $X$ the following are equivalent.
  \begin{enumerate}
  \item[a)]
    $X$ satisfies asymptotic $(\ell_p,\ell_q)$-tree estimates.
  \item[b)]
    $X$ can be embedded in a reflexive space $Z$ with an FDD $(E_i)$
    satisfying asymptotic $(\ell_p,\ell_q)$-estimates.
  \item[c)]
    $X$  is the quotient of a reflexive space $Z$ with an FDD $(E_i)$
    satisfying asymptotic  $(\ell_p,\ell_q)$-estimates.
  \item[d)]
    $X^*$ satisfies asymptotic $(\ell_{q'},\ell_{p'})$-tree
    estimates.
  \end{enumerate}
\end{cor}
\begin{remark}
  Following usual custom the range $\ell_p$, $1\leq p\leq\infty$
  really means the range $\ell_p$, $1\leq p<\infty$, $\co$.
\end{remark}
\begin{proof}
We will verify that the conditions of Theorem~\ref{T:4.6} hold for
$V=\ell_p$ and $U=\ell_q$. By the remark following Theorem~\ref{T:4.6}
these conditions are automatically satisfied if $p=\infty$ or
$q=1$. Let us now assume that $1<q\leq p<\infty$.

Since $T_{p,\gamma}$ is the $p$-convexification of $T_{1,\gamma^p}$
(see~\cite{CS}), the unit vector basis of $T_{p,\gamma}$ is
1-dominated by the unit vector basis of $\ell_p$, and hence the unit
vector basis of $\big(T_{p,\gamma}\big)^*$ 1-dominates the unit vector
basis of $\ell_{p'}$. From this one can easile deduce that every
normalized block sequence of the unit vector basis of
$\big(T_{q',\gamma}\big)^*$ 1-dominates every normalized block
sequence of the unit vector basis of $T_{p,\gamma}$.
\end{proof}

A special case of Corollary~\ref{C:4.7} solves Problem~5.4 raised
in~\cite{OS1}.
\begin{cor}\label{C:4.6}
  Let $X$ be a reflexive asymptotic $\ell_p$ space (meaning that $X$
  satisfies asymptotic $(\ell_p,\ell_p)$-tree estimates) for some $p$
  with $1\leq p\leq\infty$. Then $X$ can be embedded into a reflexive
  space with an asymptotic $\ell_p$ FDD, and $X$ is the quotient of a
  reflexive space with an asymptotic $\ell_p$ FDD.
\end{cor}


\begin{thebibliography}{MMT}
  
  
  
  
\bibitem[B]{B}  J.~Bourgain,   {\em On separable Banach spaces, universal for
  all separable reflexive spaces}, Proc. Amer. Math. Soc. {\bf 79}  (1980),
  no. 2, 241--246.
  
\bibitem[CS]{CS} P. G. Casazza and T. J. Shura
  {\em Tsirelson's Space}, Lecture Notes in Mathematics {\bf 1363},
  Springer-Verlag (1989).
  
\bibitem[FJ]{FJ} Figiel and W. B. Johnson,
  {\em A uniformly convex Banach space which contains no $\ell_p$},
  Comp. Math. {\bf 29 } (1974), 179--190
  
  
\bibitem[J]{J} W.B.~Johnson, {\em on quotients of $L_p$ which are
  quotients of $\ell_p$}, Compositio Math. {\bf 34 } (1977), 69 --
  89.

\bibitem[JZ1]{JZ1} W.B.~Johnson and M. Zippin,
  {\em On subspaces of quotients of $(\sum G_n)_{l_p}$ and $(\sum G_n)_{c_0}$},
  Proceedings of the International Symposium on Partial Differential Equations
  and the Geometry of Normed Linear Spaces, Jerusalem, 1972,
  Israel J. Math. {\bf13} (1972), 311--316 (1973).
  
  
\bibitem[JZ2]{JZ2}
  W.B. Johnson and M. Zippin,
  {\em Subspaces and quotient spaces of $(\sum G_n)_{\ell_p}$ and
    $(\sum G_n)_{c_0}$},
  Israel J. Math. {\bf 17} (1974), 50--55.
  
\bibitem[KOS]{KOS}
  H.~Knaust, E.~Odell, and Th.~Schlumprecht,
  {\em On asymptotic structure, the Szlenk index and UKK properties in
    Banach spaces},
  Positivity {\bf3} (1999), 173--199.
  
\bibitem[LT]{LT}
  J. Lindenstrauss and L. Tzafriri,
  {\em Classical Banach Spaces, I }, Springer-Verlag New York (1977).
  
\bibitem[MMT]{MMT}
  B. Maurey, V.D. Milman and N. Tomczak-Jaegermann,
  {\em Asymptotic infinite-dimensional theory of Banach spaces},
  Oper. Theory: Adv. Appl. {\bf 77} (1994), 149--175.
  
\bibitem[MT]{MT}
  V.D. Milman and N. Tomczak-Jaegermann,
  {\em Asymptotic $\ell_p$ spaces and bounded distortions},
  eds. Bor-Luh Lin and W.B. Johnson,
  Contemp. Math. {\bf 144} (1993), 173--195.
  
\bibitem[OS1]{OS1}
  E.~Odell, and Th.~Schlumprecht, {\em Trees and branches in Banach spaces},
  Trans. Amer. Math. Soc.  {\bf 354 } (2002),  no. 10, 4085--4108.
  
\bibitem[OS2]{OS2}
  E.~Odell, and Th.~Schlumprecht,  {\em A universal reflexive space
    for the class of uniformly convex Banach spaces},  to appear in
  the Mathematische Annalen.

\bibitem[P]{P}S.~Prus, {\em Finite-dimensional decompositions of Banach spaces
  with $(p,q)$-estimates},  Dissertationes Math. {\bf 263}  (1987).
  
\bibitem[Z]{Z} M. Zippin,
  {\em Banach spaces with separable duals},
  Trans. AMS, {\bf 310}, Nr. 1 (1988), 371--379.
  
  %ruffley
\end{thebibliography}
\end{document}